\documentclass[11pt]{article}
\newtheorem{lemma}{Lemma}[section]
\newtheorem{theorem}[lemma]{Theorem}
\newtheorem{proposition}[lemma]{Proposition}
\newtheorem{corollary}[lemma]{Corollary}
\newtheorem{definition}[lemma]{Definition}
\newtheorem{remark}[lemma]{Remark}
\newtheorem{remarks}[lemma]{Remarks}

\newtheorem{examples}[lemma]{Examples}

\def\square#1#2 
  {\vbox
    {\hrule
     \hbox
       {\vrule
        \vbox spread#1
          {\vfil
           \hbox spread#1{\hfil#2\hfil}
           \vfil}
        \vrule}\hrule}}
\def\sq{\square{6pt}{} }

\newenvironment{proof}{{\bf Proof}}{\hfill $ \sq $ \vskip 4mm}

\newcommand{\thlabel}[1]{\label{th:#1}}
\newcommand{\thref}[1]{Theorem~\ref{th:#1}}
\newcommand{\selabel}[1]{\label{se:#1}}
\newcommand{\seref}[1]{Section~\ref{se:#1}}

\newcommand{\prlabel}[1]{\label{pr:#1}}
\newcommand{\prref}[1]{Proposition~\ref{pr:#1}}

\newcommand{\relabel}[1]{\label{re:#1}}
\newcommand{\reref}[1]{Remark~\ref{re:#1}}

\newcommand{\delabel}[1]{\label{de:#1}}
\newcommand{\deref}[1]{Definition~\ref{de:#1}}
\newcommand{\eqlabel}[1]{\label{eq:#1}}
\newcommand{\eqref}[1]{(\ref{eq:#1})}

\newcommand{\Hom}{{\rm Hom}}
\newcommand{\HOM}{{\rm HOM}}

\newcommand{\End}{{\rm End}}

\def\lan{\langle}
\def\ran{\rangle}

\def\bu{\bullet}
\def\rhua{\hbox{$\rightharpoonup$}}
\def\lhua{\hbox{$\leftharpoonup$}}

\def\pijldubbel{\lower.2ex\vbox{\hbox{$\longrightarrow$}\vspace*{-4mm}
    \hbox{$\longrightarrow$}}}
\def\doubleright#1{{\lower.2ex\vbox{\hbox{${\smash{{\cal
M}athop{\longrightarrow}
\limits^{#1}}}$}\vspace*{-4mm}\hbox{$\longrightarrow$}}}}
\def\doublerightbis#1#2{{\lower.2ex\vbox{
\hbox{${\smash{{\cal M}athop{\longrightarrow}\limits^{#1}}}$}\vspace*{-4mm}
\hbox{${\smash{{\cal M}athop{\longrightarrow}\limits_{#2}}}$}}}}
\def\doublerightleft#1#2{{\lower.2ex\vbox{
\hbox{${\smash{{\cal M}athop{\longrightarrow}\limits^{#1}}}$}\vspace*{-4mm}
\hbox{${\smash{{\cal M}athop{\longleftarrow}\limits_{#2}}}$}}}}

\def\square#1#2 
{\vbox{\hrule\hbox{\vrule\vbox spread#1
{\vfil\hbox spread#1{\hfil#2\hfil}\vfil}\vrule}\hrule}}

\def\sq{\square{5pt}{} }

\def\ot{\otimes}

\def\doublerightleft#1#2{{\lower.2ex\vbox{
\hbox{${\smash{\mathop{\longrightarrow}\limits^{#1}}}$}\vspace*{-4mm}
\hbox{${\smash{\mathop{\longleftarrow}\limits_{#2}}}$}}}}

\input diagrams

\begin{document}
\title{Integrals, quantum Galois extensions and the
affineness criterion for quantum Yetter-Drinfel'd modules}
\author{C. Menini \\University of Ferrara\\ Departament of Mathematics \\
Via Machiavelli 35\\44100 Ferrara, Italy\\ e-mail:
men@dns.unife.it \and G. Militaru\thanks{This paper was written
while the first author was a member of G.N.S.A.G.A. with partial
financial support from M.U.R.S.T. and the second author was a
visiting professor at the University of Ferrara, supported by
C.N.R. and C.N.C.S.I.S.}
\\University of Bucharest
\\Faculty of Mathematics\\
Str. Academiei 14\\
RO-70109 Bucharest 1, Romania\\
e-mail: gmilit@al.math.unibuc.ro }

\date{}
\maketitle
\begin{abstract}
\noindent In this paper we shall generalize the notion of integral
on a Hopf algebra introduced by Sweedler, by defining the more
general concept of integral of a threetuple $(H, A, C)$, where $H$
is a Hopf algebra coacting on an algebra $A$ and acting on a
coalgebra $C$. We prove that there exists a total integral $\gamma
:C\to \Hom(C,A)$ of $(H, A, C)$ if and only if any representation
of $(H, A, C)$ is injective in a functorial way, as a
corepresentation of $C$.  In particular, the quantum integrals
associated to Yetter-Drinfel'd modules are defined. Let now $A$ be
an $H$-bicomodule algebra, ${}^H{\cal YD}_A$ the category of
quantum Yetter-Drinfel'd modules and $B=\{ a\in A \mid \sum
S^{-1}(a_{<1>}) a_{<-1>}\ot a_{<0>}= 1_H\ot a \}$,
the subalgebra of coinvariants of the Verma structure
$A\in {}^H{\cal YD}_A$. We shall prove the following affineness
criterion: if there exists $\gamma :H\to \Hom (H, A)$ a total
quantum integral and the canonical map $\beta: A\ot_B A\to H\ot A$,
$\beta (a\ot_B b)=\sum S^{-1}(b_{<1>}) b_{<-1>}\ot ab_{<0>}$
is surjective (i.e. $A/B$ is a quantum homogeneous space), then
the induction functor $-\ot_B A :{\cal M}_B \to {}^H{\cal YD}_A$
is an equivalence of categories. The affineness criteria proven by
Cline, Parshall and Scott, and independently by Oberst (for affine
algebraic groups schemes), Schneider (in the noncommutative case),
are recovered as special cases.
\end{abstract}

\section{Introduction}\selabel{0}
The integrals for Hopf algebras were introduced in two fundamental
papers: by Larson and Sweedler in \cite{LS} for the finite case,
and by Sweedler in \cite{Sw} for the infinite case. Initially
introduced in order to generalize the Haar measure on a compact
group, the integrals have proven to be a powerful instrument in
the classic theory of Hopf algebras, beginning with representation
theory and ending with the classification theory for finite
dimensional Hopf algebras. Recently, arising from an idea of
Drinfel'd, the integrals were introduced for Hopf algebras in
various braided categories: abelian and rigid (\cite{Ly}) or rigid
with split idempotents (\cite{BKLT}). At this level (mutatis
mutandis, the definition is fundamentally the one given by
Sweedler), integrals have proven to be essential instruments in
constructing invariants of surgically presented 3-manifolds or
3-dimensional topological quantum field theories (\cite{Kerler},
\cite{Ku}, \cite{T}).

In the first part of this paper we shall introduce the more general concept
of integral associated to a threetuple $(H, A, C)$ called Doi-Koppinen
datum,
consisting of a Hopf algebra $H$ which coacts on an algebra
$A$ and acts on a coalgebra $C$. As a major application, the
quantum integrals associated to Yetter-Drinfel'd modules
${}^H{\cal YD}_H$, are introduced.
The transition from the classic integrals of Sweedler, which are elements
$\varphi \in H^*$ invariant to convolution
(or equivalently $H$-colinear maps $\varphi: H\to k$),
to the quantum integrals, which are maps $\gamma :H\to \End(H)$
satisfying the condition
$$
\sum g_{(1)}\ot \gamma(g_{(2)})(h)=
\sum S^{-1} \Bigl (\{\gamma (g)(h_{(1)})\}_{(3)}\Bigl)
h_{(2)}\{\gamma(g)(h_{(1)})\}_{(1)}\ot \{\gamma(g)(h_{(1)})\}_{(2)}
$$
for all $g$, $h\in H$, is a long way which needs to be explained.

First of all, the integrals on a Hopf algebra $H$ and the more
general ones introduced by Doi (\cite{D0}) for an $H$-comodule
algebra $A$, have strong ties to ${\cal M}^H$, the
corepresentations of $H$, and to the representations of the pair
$(H, A)$, being the category of relative Hopf modules ${\cal
M}^H_A$. The starting point for this paper is presented in
\seref{1.2}: the existence of an integral in the sense of Doi
(classic, if we consider $A=k$) is the necessary and sufficient
criterion for the existence of a natural transformation between
two functors linking ${\cal M}^H_A$ to ${\cal M}^H$ (see
\thref{1.2.12}).

This categorical point of view towards integrals will allow us to
correctly define the integrals associated to a
Doi-Koppinen datum $(H, A, C)$.
We shall thus arrive at the end of a road initiated in \cite{CMZ2},
where the $H$-integrals of a Doi-Koppinen datum were introduced
in relation to Frobenius type theorems, and continued in \cite{CMZ3},
\cite{CIMZ}, where $A$-integrals were defined as tools for proving
separability theorems.
Let ${}^C{\cal M}(H)_A$ be the
category of representations of $(H, A, C)$: an object of it
(also called a Doi-Koppinen module) is a $k$-module with an
$A$-action and a compatible $C$-coaction.
In view of the above, an integral of
$(H, A, C)$ will be a map, the existence of which is
equivalent with the existence of a natural
transformation between the functors
$F_A\circ (C \ot \bullet)\circ F^C$ and
$F_A\circ 1_{{}^C{\cal M}(H)_A}: {}^C{\cal M}(H)_A\to
{}^C{\cal M}$,
where $F_A:{}^C{\cal M}(H)_A \to {}^C{\cal M}$ and
$F^C:{}^C{\cal M}(H)_A\to {\cal M}_A$
are the corresponding forgetful functors.

The definition of the integral of $(H,A,C)$, i.e. a map $\gamma
:C\to \Hom(C,A)$ satisfying the equation \eqref{TI1}, is given in
\deref{2.1.4a}, and its characterization in \thref{doi}, which can
be interpreted as follows: there exists a total integral $\gamma
:C\to \Hom(C,A)$ of $(H, A, C)$ if and only if any Doi-Koppinen
module is relative injective (injective, if we work over a field)
in a functorial way, as a left $C$-comodule. A key result is
\thref{mm3} : if there exists a total integral $\gamma :C\to
\Hom(C,A)$ of $(H,A,C)$, then $C\ot A$ is a generator in the
category ${}^C{\cal M}(H)_A$. As explained in previous
publications (\cite{CIMZ}, \cite{CMZ1}), ${}^C{\cal M}(H)_A$
unifies modules, comodules, Sweedler's Hopf modules, relative Hopf
modules, graded modules, Long dimodules and Yetter-Drinfel'd
modules. In particular, by applying the above results we obtain
the definition and the characterization theorem for quantum
integrals, being the integrals which correspond to the
Yetter-Drinfel'd modules ${}^H{\cal YD}_H$: If there exists a
total quantum integral $\gamma :H\to \End(H)$, then any
Yetter-Drinfel'd module $M\in {}^H{\cal YD}_H$ is relative
injective as an $H$-comodule. In particular, if $H$ is finite
dimensional over a field $k$, there exists a total quantum
integral $\gamma :H\to \End(H)$ if and only if any representation
of the Drinfel'd double $D(H)$ is injective in a functorial way,
as an $H^*$-module.

In the second part of the paper we introduce the notion of quantum
Galois extensions and we prove a criterion for affineness in a
quantum version. Let us explain the terminology and
justify the usage of the term "quantum". We shall begin by
recalling the following powerful theorem given by Schneider (\cite{Sch})
\footnote{For $H$ and $A$ commutative, this result was proved before
by Doi in \cite[Theorem 3.2]{D0}.} (presenting an equivalent
right-left version of it):

\begin{theorem}\thlabel{mmsch}
Let $H$ be a Hopf algebra with a bijective antipode over a field $k$,
$A$ a left $H$-comodule algebra and $B=A^{{\rm co}(H)}$ its subalgebra
of coinvariants. The following statements are equivalent
\begin{enumerate}
\item
\begin{enumerate}
\item there exists $\varphi :H\to A$ a total integral;
\item the canonical map ${\rm can} :A\ot_B A \to H\ot A$,
${\rm can}(a\ot_B b)=\sum b_{<-1>}\ot ab_{<0>}$, is surjective;
\end{enumerate}
\item the induction functor
$-\ot_B A :{\cal M}_B \to {}^H{\cal M}_A$
is an equivalence of categories;
\item
\begin{enumerate}
\item $A$ is faithfully flat as a left $B$-module;
\item $A/B$ is an $H$-Galois extension, i.e. ${\rm can}$ is bijective.
\end{enumerate}
\end{enumerate}
\end{theorem}
A few comments on this result:\\
$3)\Rightarrow 1)$
was proven for the case when $k$ is a field,
\footnote{Recently, in \cite{MSTW}, it was proven using other
techniques that this implication holds for commutative QF-rings.}
based on a result of Takeuchi \cite{Take}: over a field,
an $H$-comodule is injective if and only if it is coflat.
All the other implications hold true over a
commutative ring, if we additionally assume that $H$ is projective
over $k$.
The equivalence $2) \Leftrightarrow 3)$ is a general imprimitivity
theorem and its proof is standard from the point of view
of category theory: a pair of adjoint functors (as
$-\ot_B A$ and $(-)^{{\rm co}(H)}$ are)
gives an equivalence of categories iff one of them is faithfully
exact (or both of them are exact) and the adjunction maps in
the key objects of categories ($B$ in ${\cal M}_B$ and
$H\ot A$ in ${}^H{\cal M}_A$) are bijective.
\footnote{For generalizations of this equivalence we refer
to \cite[Theorem 5.6]{Brzezinski00}, \cite[Theorem 3.10]{Br}
or \cite[Theorem 2.8]{CR}.}
The main part of the theorem is $1) \Rightarrow 2) $, which is the
non-commutative version of the affineness criterion
for affine algebraic groups schemes given before by
Cline, Parshall and Scott \cite{CPS} and independently by Oberst
\cite{Ob}.

Our intention is to "quantize" this result. The first step to be
taken is to quantize the category ${}^H{\cal M}_A$ in a coherent
way. In Hopf algebras theory, quantization means (roughly
speaking) a deformation of the enveloping algebra of a semisimple
Lie algebra ${\ell}$ using a parameter $q$, in order to obtain a
noncommutative noncocommutative Hopf algebra (a quantum group)
$U_q({\ell})$ (\cite{Dr}, \cite{J}). The results obtained for the
new object $U_q({\ell})$ in the framework of representation theory
(\cite{Lu1}, \cite{Rosso}) will generalize the results from the
classic case $U({\ell})$. In order to quantize the category
${}^H{\cal M}_A$, the part of the parameter $q$ mentioned above
will be played by a new coaction (i.e. a family of parameters,
possible infinite, if $k$ is a field), this time to the right,
$\rho^r :A\to A\ot H$. Thus $A$ will be a $H$-bicomodule algebra
and the category of representations will be denoted by ${}^H{\cal
YD}_A$, the category of quantum Yetter-Drinfel'd modules
introduced in \cite{CMZ1}. If $A=H$ and $\rho^l=\rho^r=\Delta$
then ${}^H{\cal YD}_H$ is the category of Yetter-Drinfel'd (or
crossed) modules introduced in \cite{Dr}, \cite{Y}: if $H$ is
finite dimensional then ${}^H{\cal YD}_H\cong {\cal M}_{D(H)}$,
where $D(H)$ is the Drinfel'd double of $H$ (\cite{Majid91}). The
category ${}^H{\cal YD}_H$ plays an important role in the quantum
Yang-Baxter equation, low dimensional topology or knot theory
(\cite{Kassel}, \cite{Majid}). In \cite{EK}, the category of
Yetter-Drinfel'd modules\footnote{The authors consider the
left-left equivalent version of ${}^H{\cal YD}_H$ and the object
of it are called $H$-dimodules.} was used as the fundamental tool
in the construction of the {\sl dequantization} functor {\sl DQ},
which is an equivalence of categories from the category of
quantized universal enveloping algebras to the category of Lie
bialgebras over $k[[h]]$.\\ On the other hand, if $\rho_r :A\to
A\ot H$ is the trivial coaction, that is $\rho^r(a)=a\ot 1_H$,
then ${}^H{\cal YD}_A={}^H{\cal M}_A$, the category of classical
relative Hopf modules.

Hence, ${}^H{\cal YD}_A$ is a category containing the category of
Yetter-Drinfel'd modules ${}^H{\cal YD}_H$ as a particular case
and, on the other hand, ${}^H{\cal M}_A$ is obtained from
${}^H{\cal YD}_A$ by trivializing the right coaction of $H$ on
$A$. We can therefore view the category ${}^H{\cal YD}_A$ as a
quantization of the category of relative Hopf modules ${}^H{\cal
M}_A$. For this reason, throughout this paper we shall call the
objects of ${}^H{\cal YD}_H$ Yetter-Drinfel'd modules, while the
objects of the more general category ${}^H{\cal YD}_A$ shall be
called quantum Yetter-Drinfel'd modules.

Now we can extend \thref{mmsch} from relative Hopf modules
to quantum Yetter-Drinfel'd modules: this is what we shall do in
\seref{3}. With the exception of $3) \Rightarrow 1)$
(which remains an open problem), all the other implications
maintain their validity for the category ${}^H{\cal YD}_A$.
Let
$$
B=\{ a\in A \mid \sum S^{-1}(a_{<1>}) a_{<-1>}\ot a_{<0>}=
1_H\ot a \},
$$
be the subalgebra of quantum coinvariants, which are the
coinvariants of the Verma structure
$A\in {}^H{\cal YD}_A$.
\prref{gal70} shows that the induction functor
$-\ot_B A :{\cal M}_B \to {}^H{\cal YD}_A$
is an equivalence of categories if and only if $A/B$ is a
faithfully flat quantum Galois extension, i.e. the
canonical map $\beta: A\ot_B A\to H\ot A$,
$\beta (a\ot_B b)=\sum S^{-1}(b_{<1>}) b_{<-1>}\ot ab_{<0>}$ is
bijective. \thref{mm80} proves the quantum affineness criterion:
if there exists $\gamma :H\to \Hom (H, A)$ a total quantum
integral and
the canonical map $\beta: A\ot_B A\to H\ot A$ is surjective,
then the induction functor
$-\ot_B A :{\cal M}_B \to {}^H{\cal YD}_A$
is an equivalence of categories.

\section{Preliminary results}\selabel{1}
Throughout this paper, $k$ will be a commutative ring with unit. Unless
specified otherwise, all modules, algebras, coalgebras, bialgebras,
tensor products and homomorphisms are over $k$. For a $k$-algebra $A$,
${{\cal M}}_A$ (resp. ${}_A{\cal M}$) will be the category of right (resp.
left) $A$-modules and $A$-linear maps. $H$ will be a Hopf
algebra over $k$, and we will use Sweedler's
sigma-notation extensively. For example,
if ($C,\Delta $) is a coalgebra,
then for all $c\in C$ we write
$$
\Delta (c)=\sum c_{(1)}\ot c_{(2)}\in C\ot C, \quad
(\Delta \ot {\rm Id})\Delta (c)= ({\rm Id} \ot \Delta)\Delta (c)=
\sum c_{(1)}\ot c_{(2)}\ot c_{(3)}.
$$
If ($M,\rho_M$) is a left $C$-comodule, then we write
$$
\rho_M(m)=\sum m_{<-1>}\otimes m_{<0>} \in C\ot M,$$
and
$$
(\Delta \ot {\rm Id})\rho_M(m)= ({\rm Id} \ot \rho_M)\rho_M (m)=
\sum m_{<-2>}\ot m_{<-1>}\ot m_{<0>} \in C\ot C\ot M
$$
for $m\in M$.
${}^C{{\cal M}}$ (resp. ${\cal M}^C$)
will be the category of left (resp. right) $C$-comodules and
$C$-colinear maps. If $M$ is a left $C$-comodule then $C\ot M$ is
also a left $C$-comodule via $\rho_{C\ot M}:= \Delta \ot Id_M$ and
$\rho_M :M\to C\ot M$ is a left $C$-colinear map. A left $C$-comodule
$M$ is called {\it relative injective}  if for any
$k$-split monomorphism $i:U\to V$ in ${}^C{{\cal M}}$ and for any
$C$-colinear map $f:U\to M$, there exists a $C$-colinear map
$g:V\to M$ such that $g\circ i=f$. This is equivalent (\cite{D0})
to the fact that $\rho_M:M\to C\ot M$ splits in ${}^C{{\cal M}}$,
i.e. there exists
a $C$-colinear map $\lambda_M :C\ot M\to M$ such that
$\lambda_M \circ \rho_M= {\rm Id}$. Of course, if $k$ is a field, $M$ is
relative injective if and only if it is an injective object
in ${}^C{{\cal M}}$.\\
The dual $C^*=\Hom(C,k)$ of a $k$-coalgebra $C$ is a $k$-algebra. The
multiplication on $C^*$ is given by the convolution
$$\lan f*g,c\ran=\sum \lan f,c_{(1)}\ran\lan g,c_{(2)}\ran,$$
for all $f,g\in C^*$ and $c\in C$. $C$ is a $C^*$-bimodule: the
left and right action are given by the formulas
\begin{equation}\eqlabel{1.0a}
c^*\rhua c=\sum \lan c^*,c_{(2)}\ran c_{(1)}~~{\rm and}~~
c\lhua c^*= \sum \lan c^*,c_{(1)}\ran c_{(2)}
\end{equation}
for $c^*\in C^*$ and $c\in C$. This also holds for $C$-comodules:
for example, if ($M,\rho_M)$ is a left $C$-comodule, then it
becomes a right $C^*$-module by $$m\cdot c^*=\sum\lan c^*,
m_{<-1>}\ran m_{<0>},$$ for all $m\in M$ and $c^*\in C^*$.\\ An
algebra $A$ that is also a left $H$-comodule is called a left
$H$-comodule algebra if the comodule structure map $\rho_A :A\to
H\ot A$ is an algebra map. This means that $$\rho_A(ab)=\sum
a_{<-1>}b_{<-1>}\ot a_{<0>}b_{<0>}~~{\rm and}~~ \rho_A(1_A)=1_H\ot
1_A$$ for all $a,b\in A$. This is equivalent to the fact that $A$
is an algebra in the monoidal category ${}^H{\cal M}$ of left
$H$-comodules. \\ Similarly, a coalgebra that is also a right
$H$-module is called a right $H$-module coalgebra if $C$ is a
coalgebra in the monoidal category ${\cal M}_H$ of right
$H$-modules, or equivalent $$\Delta_C(c\cdot h)=\sum c_{(1)}\cdot
h_{(1)}\ot c_{(2)}\cdot h_{(2)} ~~{\rm and}~~ \varepsilon_C(c\cdot
h)=\varepsilon_C(c)\varepsilon_H(h),$$ for all $c\in C$, $h\in H$.

\subsection{Doi-Koppinen modules: functors and structures}\selabel{1.1}
Let $H$ be a Hopf algebra, $A$ a left $H$-comodule algebra and
$C$ a right $H$-module coalgebra.
The threetuple ${\cal O}=(H,A,C)$ is called a {\it Doi-Koppinen datum}.
In order to study these general objects, we
have to define their representations: a representation of
${\cal O}=(H,A,C)$ (called also a right-left {\it Doi-Koppinen module})
is a
$k$-module $M$ that has the structure of right $A$-module and
left $C$-comodule, such that the
following compatibility relation holds
\begin{equation}\eqlabel{comp}\eqlabel{1.1a}
\rho_M(m a)=\sum m_{ <-1>}\cdot a_{<-1>} \otimes m_{<0>}a_{<0>},
\end{equation}
for all $a\in A$, $m\in M$. ${}^C{\cal M}(H)_A$ will be the
abelian category of right-left Doi-Koppinen modules and
$A$-linear, $C$-colinear maps, as it was introduced by Y. Doi in
\cite{D2} and independently by M. Koppinen in \cite{K}. Let $F^C:
{}^C{\cal M}(H)_A\to {\cal M}_A$ be the forgetful functor which
forgets the $C$-coaction and $$C\ot \bullet :\ {\cal M}_A\to
{}^C{\cal M}(H)_A, \quad M\to C\ot M$$ its right adjoint, where
the structure maps on $C\ot M$ are given by
\begin{eqnarray}
(c\ot m)\cdot a&=&\sum c\cdot a_{<-1>}\ot ma_{<0>}\eqlabel{1.1}\\
\rho_{C\ot M}(c\otimes m)&=&\sum c_{(1)}\otimes c_{(2)}\otimes m \eqlabel{1.2}
\end{eqnarray}
for any $c\in C,~a\in A$ and $m\in M$.
The unit of the adjoint pair $(F^C, C\ot \bullet)$ is precisely
(\cite{CIMZ}, \cite{Mi})
$$\rho:\ 1_{{}^C{\cal M}(H)_A} \rightarrow (C\ot \bullet )\circ F^C, $$
the $C$-coaction $\rho_M:\ M\to C\ot M$ on any Doi-Koppinen module $M$;
therefore $\rho_M$ is $A$-linear and $C$-colinear and
can be viewed as a natural transformation between the
functors $1_{{}^C{\cal M}(H)_A}$ and
$(C\ot \bullet )\circ F^C$. $A$ is a right $A$-module, so
$C\ot A$ is a Doi-Koppinen module via:
\begin{eqnarray}
(c\ot b)a&=& \sum ca_{<-1>}\ot ba_{<0>}\eqlabel{1.3}\\
\rho^{l}_{C\ot A}(c\ot b)&=& \sum c_{(1)}\ot c_{(2)}\ot b\eqlabel{1.4}
\end{eqnarray}
Furthermore, $C\ot A$ is also a right $C$-comodule via
\begin{equation}\eqlabel{sase}
\rho^{r}_{C\ot A}: C\ot A \to C\ot A\ot C, \quad
\rho^{r}_{C\ot A}(c\ot a)=\sum c_{(1)}\ot a_{<0>}\ot c_{(2)}S(a_{<-1>})
\end{equation}
for all $a\in A$, $c\in C$.
$C\ot C\ot A= \Bigl( (C\ot\bullet)\circ F^C \Bigl) (C\ot A)$ is
an object in ${}^C{\cal M}(H)_A$ and also has a right $C$-comodule
structure via:
\begin{eqnarray}
(c\ot d\ot b)a&=& \sum ca_{<-2>}\ot da_{<-1>}\ot ba_{<0>}\eqlabel{1.5}\\
\rho^{l}_{C\ot C\ot A}(c\ot d\ot b)&=&
\sum c_{(1)}\ot c_{(2)}\ot d\ot b\eqlabel{1.6}\\
\rho^{r}_{C\ot C\ot A}(c\ot d\ot b)&=&c\ot \rho^{r}_{C\ot A} (d\ot b)=
\sum c\ot d_{(1)}\ot b_{<0>}\ot d_{(2)}S(b_{<-1>}) \eqlabel{11}
\end{eqnarray}
Now, let $F_A :{}^C{\cal M}(H)_A\to {}^C{\cal M}$ be the other
forgetful functor, which forgets the $A$-action and
$$ \bullet \ot A :{}^C{\cal M}\to {}^C{\cal M}(H)_A, \quad
N\to N\ot A$$
its left adjoint, where for $N\in {}^C{\cal M}$,
$N\ot A \in {}^C{\cal M}(H)_A$ via the structures
\begin{eqnarray}
(n\ot a)\cdot b&=&n\ot ab\\
\rho_{N\ot A}(n\otimes a)&=&
\sum n_{<-1>}a_{<-1>}\ot n_{<0>}\ot a_{<0>}  \eqlabel{opt}
\end{eqnarray}
for any $a, ~b\in A$ and $n\in N$. $C$ is a left $C$-comodule via
$\Delta$; hence $C\ot A$ can be also viewed as a Doi-Koppinen module via
\begin{eqnarray}
(c\ot b)\cdot^{\prime}a&=& c\ot ba\eqlabel{1.37}\\
\rho^{\prime l}_{C\ot A}(c\ot b)&=&
\sum c_{(1)}\cdot b_{<-1>}\ot c_{(2)}\ot b_{<0>}\eqlabel{1.48}
\end{eqnarray}
for all $c\in C$, $a$, $b\in A$. These two types of Doi-Koppinen module
structures on $C\ot A$, coming from \eqref{1.37} and \eqref{1.48}
or from \eqref{1.3} and \eqref{1.4}, are isomorphic:
more precisely, the map
\begin{equation}\eqlabel{mm1.5}
u: C\ot A \to C\ot A, \quad u(c\ot a)= \sum c a_{<-1>}\ot a_{<0>}
\end{equation}
is an isomorphism of Doi-Koppinen modules (\cite{CR}) with an inverse
given by
$$
u^{-1}: C\ot A \to C\ot A, \quad u^{-1}(c\ot a)= \sum c S(a_{<-1>})\ot a_{<0>}
$$
The algebra $C^*$ is a left $H$-module algebra; the $H$-action
is given by the formula
$$\lan h\cdot c^*,c\ran=\lan c^*,c\cdot h\ran$$
for all $h\in H$, $c\in C$ and $c^*\in C^*$. The smash product $A\#C^*$ is
equal to $A\ot C^*$ as a $k$-module, with the multiplication defined by
\begin{equation}\eqlabel{smash}
(a\#c^*)(b\#d^*)=\sum a_{<0>}b\#c^{*}* (a_{<-1>}\cdot d^*),
\end{equation}
for all $a, b\in A$, $c^*, d^*\in D^*$. Recall that we have a
natural functor
$P:\ {}^C{\cal M}(H)_A\to {\cal M}_{A\# C^*}$
sending a Doi-Koppinen module $M$
to itself,
with the right $A\# C^*$-action given by
\begin{equation}\eqlabel{smash2}
m\cdot (a\# c^*)=\sum\lan c^*,m_{<-1>}\ran m_{<0>}a
\end{equation}
for any $m\in M$, $a\in A$ and $c^*\in C^*$. $P$ is an equivalence of
categories if $C$ is finitely generated and projective as a $k$-module
(\cite{D2}).\\
As in \cite{CIMZ}, the following right $C^*$-module on $\Hom (C, A)$
will play a key role
\begin{eqnarray}
(f\cdot c^*)(c)&=&\sum f\bigl(c_{(1)}\bigr)_{<0>}\lan c^*,
c_{(2)}\cdot f\bigl(c_{(1)}\bigr)_{<-1>}\ran\eqlabel{1.10c}
\end{eqnarray}
for any $f\in \Hom(C,A),~c^*\in C^* ~{\rm and}~ c\in C$.
The above
right $C^*$-action on $\Hom(C,A)$ is very natural:
$\Hom(C,A)$ has an algebra structure which is
the right-left version of the smash product in the sense of
Koppinen (\cite{K}), i.e. the multiplication is given by
\begin{equation}\eqlabel{Ko}
(f\bu g)(c)=\sum f\left(c_{(1)}\right)_{<0>}
g\left( c_{(2)}\cdot f\left( c_{(1)}\right)_{<-1>}\right)
\end{equation}
for all $f$, $g\in \Hom(C,A)$, $c\in C$. Moreover, if $C$ is
finitely generated and projective, the canonical map
$i:\ A\#C^*\to \Hom(C,A)$ given by
$i(a\# c^*)(c)=\lan c^*,c \ran a$ is an algebra isomorphism.
Now, $\Hom(C,A)$ contains $C^*$ as a subalgebra via
$j : C^* \to \Hom (C, A)$, $j(c^*)(c):=\lan c^*,c \ran $,
and the
$C^*$-action given by \eqref{1.10c}
is exactly the structure induced by the usual restriction of scalars
via $j: C^*\to \Hom(C,A)$. \\
Let $A$ be a left $H$-comodule algebra and $C=H$, viewed as a right
$H$-module coalgebra via the multiplication on $H$. Then
${}^H{\cal M}(H)_A={}^H{\cal M}_A$, the category of (right-left)
relative Hopf-modules. We can also define ${}^H_A{\cal M}$, the
category of (left-left) relative Hopf modules: an object in this
category is a $k$-module $M$ which has a left $A$-module structure
and a left $H$-comodule structure such that the
following compatibility relation holds
\begin{equation}\eqlabel{compll}\eqlabel{1.1ab}
\rho_M(a m)=\sum a_{ <-1>}\cdot m_{<-1>} \otimes a_{<0>}m_{<0>},
\end{equation}
for all $a\in A$, $m\in M$. If $S$ is bijective, then $H^{{\rm op}}$
is a Hopf
algebra with $S^{-1}$ as an antipode, and there exists an equivalence of
categories
${}^H{\cal M}_A\cong {}^{H^{{\rm op}}}_{A^{{\rm op}}}{\cal M}$. \\
Similarly, if $A$ is a right $H$-comodule algebra we can define the
categories ${\cal M}^H_A$ and ${}_A {\cal M}^H$. For instance,
an object in ${\cal M}^H_A$ is a right $A$-module and right
$H$-comodule $M$ such that the
following compatibility relation holds
\begin{equation}\eqlabel{comprr}\eqlabel{1.1ac}
\rho_M(m a)=\sum m_{ <0>}\cdot a_{<0>} \otimes m_{<1>}a_{<1>},
\end{equation}
for all $a\in A$, $m\in M$. The category ${\cal M}^H_A$, together
with all its left-right equivalent versions, was introduced in
\cite{Doi1983}, and proved to be a unifying framework, at the
level of Hopf algebras, for problems arising from different and
apparently unrelated fields like: affine algebraic groups or more
generally affine scheme, Lie algebras, compact topological groups,
groups representation theory, Galois theory, Clifford theory or
graded rings theory (see, for instance, \cite{Sch}, \cite{Sch1}).

\subsection{Natural transformations versus integrals}\selabel{1.2}
In this section we shall present a point of view which is
essential for the rest of the paper: specifically, we shall prove
that the existence of an integral on a Hopf algebra is a necessary
and sufficient criterion for constructing a natural transformation
between two functors. This observation will allow us to define the
general concept of an integral associate to a Doi-Koppinen datum
$(H, A, C)$.\\ Let $H$ be a Hopf algebra over a field $k$. We
recall (\cite{Sw}) that a right integral on $H$ is an element
$\varphi \in H^*$ such that $\varphi h^*= \lan h^*, 1_H\ran
\varphi$ for all $h^*\in H^*$ . This is equivalent to the fact
that $\varphi: H\to k$ is right $H$-comodule map, where $k$ has
the trivial right $H$-comodule structure. If a right (or left)
integral exists, then the antipode of $H$ is bijective
(\cite{Radford}). \\ Doi (\cite{D0}) generalizes this concept in
the obvious way as follows: let $A$ be a right $H$-comodule
algebra. A map $\varphi :H\to A$ is called an {\it integral}
(\cite{D0}) if $\varphi$ is right $H$-colinear. Furthermore,
$\varphi$ is called a {\it total integral} if additionally
$\varphi(1_H)=1_A$. The criterion for the existence of a total
integral is given by Theorem 1.6 of \cite{D0} (we shall present
only its essential part):

\begin{theorem}\thlabel{1.2.11}
Let $A$ be a right $H$-comodule algebra. The following are equivalent:
\begin{enumerate}
\item there exists a total integral $\varphi :H\to A$;
\item any Hopf module $M\in {\cal M}_A^H$ is relative injective as
a right $H$-comodule, i.e. the
right $H$-coaction $\rho_M :M\to M\ot H$ splits in the category
${\cal M}^H$ of right $H$-comodules;
\item $\rho_A :A\to A\ot H$ splits in the category
${\cal M}^H$ of right $H$-comodules.
\end{enumerate}
\end{theorem}

We will explain now what is behind the proof of this
characterization theorem. As $(A, \rho_A)\in {\cal M}_A^H$, $2)
\Rightarrow 3) $ is trivial. Now, if $\lambda_A$ is an
$H$-colinear retraction of $\rho_A$, then $\varphi :H\to A$,
$\varphi (h):=\lambda_A (1_A\ot h)$ is a total integral, i.e $3)
\Rightarrow 1) $ follows. The implication $1) \Rightarrow 2) $ is
also easy: if $\varphi$ is a total integral then
$$\lambda_M^{\varphi}: M\ot H\to M, \quad \lambda_M^{\varphi}
(m\ot h)=\sum m_{<0>}\varphi (S(m_{<1>}) h)$$ is an $H$-colinear
retraction of $\rho_M$.\\ There is however more to be read between
the lines of this proof, and this is a main starting point for
this paper. The character of $\lambda$ is functorial: more
precisely, if $f:M\to N$ is a morphism in ${\cal M}_A^H$ (i.e. $f$
is $A$-linear and $H$-colinear), then the diagram $$
\begin{diagram}
M\ot H&\rTo{\lambda_M^{\varphi}}&M\\
\dTo{f\ot {\rm Id}}&&\dTo_{f}\\
N\ot H&\rTo{\lambda_N^{\varphi}}&N\\
\end{diagram}
$$
is commutative. Hence, $\lambda^{\varphi}$ is a natural transformation
$$
\lambda^{\varphi}:F_A\circ (\bullet \ot H)\circ F^H \to
F_A\circ 1_{{\cal M}_A^H}
$$
where $F_A:{\cal M}_A^H \to {\cal M}^H$
(resp. $F^H:{\cal M}_A^H\to {\cal M}_A$) are the corresponding
forgetful functors. Now we view the right $H$-coaction
$\rho$ as a natural transformation
$$
\rho:F_A\circ 1_{{\cal M}_A^H} \to
F_A\circ (\bullet \ot H)\circ F^H.
$$
Bearing in mind the above, the theorem of Doi can be restated as follows:

\begin{theorem}\thlabel{1.2.12}
Let $A$ be a right $H$-comodule algebra. The following are equivalent:
\begin{enumerate}
\item there exists a total integral $\varphi :H\to A$;
\item there exists a natural transformation
$\lambda:F_A\circ (\bullet \ot H)\circ F^H \to
F_A\circ 1_{{\cal M}_A^H}$
that splits
$\rho:F_A\circ 1_{{\cal M}_A^H} \to
F_A\circ (\bullet \ot H)\circ F^H$;
\item $\rho_A :A\to A\ot H$ splits in the category
${\cal M}^H$ of right $H$-comodules.
\end{enumerate}
\end{theorem}

\begin{remarks} \label{zof}
\rm 1.
The above theorem is still valid leaving aside the normalizing
condition $\varphi(1_H)=1_A$. More exactly, there exists an
integral $\varphi :H\to A$ if and only if there exists
$\lambda:F_A\circ (\bullet \ot H)\circ F^H \to
F_A\circ 1_{{\cal M}_A^H}$
a natural transformation. In particular, in the classic case
corresponding to $A=k$, we obtain that there exists a right
integral $\varphi :H\to k$ on $H$ if and only if
there exists a natural transformation
$\lambda:(\bullet \ot H)\circ F^H \to 1_{{\cal M}^H}$.
Furthermore, $\varphi (1_H)=1$ if and only if $\lambda$ splits
$\rho: 1_{{\cal M}^H} \to (\bullet \ot H)\circ F^H$. This
is equivalent (\cite{CIMZ}, \cite{Mi}) to the fact that
the forgetful functor
$F^H :{\cal M}^H\to {\cal M}_k$ is separable, which is another way
of formulating Maschke's theorem for Hopf algebras (\cite{Sw}). \\
2. Let $A$ be a left $H$-comodule algebra. The version of
\thref{1.2.12} for the category ${}^H_A{\cal M}$ is still true.
In this case the $H$-colinear split of $\rho_M :M\to H\ot M$ associated
to a left total integral $\varphi :H\to A$ is given by the formula:
\begin{equation}\eqlabel{1ss}
\lambda_M^{\prime}: H\ot M\to M, \quad
\lambda_M^{\prime}(h\ot m)=\sum \varphi (h S(m_{<-1>}))m_{<0>}
\end{equation}
for all $h\in H$, $m\in M$. \\
Now, if we deal with ${}^H{\cal M}_A$ we have to assume that the
antipode of $H$ is bijective, in order to be able to construct a
splitting for $\rho_M$. In this case, the only possible way of
constructing
a left $H$-colinear split of $\rho_M :M\to H\ot M$ seems to be the one
given by the formula:
\begin{equation}\eqlabel{1ds}
\lambda_M^{\prime\prime}: H\ot M\to M, \quad
\lambda_M^{\prime\prime}(h\ot m)=\sum m_{<0>}\varphi (S^{-1}(m_{<-1>})h)
\end{equation}
for all $h\in H$, $m\in M$, where $S^{-1}$ is the inverse of $S$.
Of course, in the trivial case $A=k$ which corresponds to classic
integrals, this is not really a restriction, as the existence of
a left integral on $H$ ensures the bijectivity of
the antipode (\cite{Radford}). We shall see however that,
even in the case
${}^H{\cal M}_A$, the restriction "$S$ bijective" can be left
behind (moreover, we can do the same with the condition that an
antipode exists), if we replace the Doi integrals $\varphi: H\to A$
with maps $\gamma :H\to \Hom(H, A)$. For the latter, the
split of $\rho_M$ is given by
$$
\lambda_M^{\prime\prime}(h\ot m)=\sum m_{<0>}\gamma (h)(m_{<-1>})
$$
for all $h\in H$, $m\in M$.
\end{remarks}

\subsection{Quantum Yetter-Drinfel'd modules}\selabel{1.4}
Let $H$ be a Hopf algebra with a bijective antipode,
$A$ an $H$-bicomodule algebra, and $C$
an $H$-bimodule coalgebra, this means that $A$ is an algebra in
the monoidal category ${}^H{\cal M}^H$ of $H$-bicomodules and
$C$ is a coalgebra in the category ${}_H{\cal M}_H$ of
$H$-bimodules. The left and right $H$-coaction
on $A$ are denoted by $\rho^l :A\to H\ot A$,
$\rho^l (a)=\sum a_{<-1>}\ot a_{<0>}$ and
$\rho^r :A\to A\ot H$, $\rho^r (a)=\sum a_{<0>}\ot a_{<1>}$.
The threetuple $G=(H,A,C)$ was called in \cite{CMZ1} a
{\em Yetter-Drinfel'd datum}.
Let $G=(H,A,C)$ be a Yetter-Drinfel'd datum. A representation
of $G$, also called a {\em crossed $G$-module} or a
{\em quantum Yetter-Drinfel'd module}, is a  $k$-module $M$ that
is at the same time a right $A$-module and a left $C$-comodule
such that
\begin{equation}\eqlabel{gal1131}
\sum m_{<-1>}a_{<-1>}\ot m_{<0>}\cdot a_{<0>}=
\sum a_{<1>} (m\cdot a_{<0>})_{<-1>}\ot (m\cdot a_{<0>})_{<0>}
\end{equation}
or, equivalently,
\begin{equation}\eqlabel{:2}
\rho_M(ma)=\sum S^{-1}(a_{<1>})\cdot m_{<-1>}\cdot
a_{<-1>}\otimes m_{<0>}a_{<0>}
\end{equation}
for all $m\in M$ and $a\in A$. The category of (right-left)
crossed $G$-modules and $A$-linear, $C$-colinear maps will be
denoted by ${}^C{\cal YD}(H)_A$ and was introduced in \cite{CMZ1}.
It follows easily that, for $C=A=H$, we obtain the classical
Yetter-Drinfel'd modules ${}^H{\cal YD}_H$ (\cite{Y}, \cite{RT}).
\\ In a similar way, we can introduce left-right, right-right and
left-left crossed $G$-modules. The corresponding categories are
${}_A{\cal YD}(H)^C$, ${\cal YD}(H)_A^C$ and ${}^C_A{\cal YD}(H)$.
There exist relationships (we refer to \cite{CMZ1} for full
details) between the four types of crossed $G$-modules, given by
the following equivalence of categories $$ {}_A{\cal
YD}(H)^C\cong\; {}_{A}^{C^{{\rm cop}}}{\cal YD}(H)\cong\; {\cal
YD}(H^{{\rm op\; cop}})^{C}_{A^{{\rm op}}}\cong \; {}^{C^{{\rm
cop}}}{\cal YD}(H^{{\rm op\; cop}})_{A^{{\rm op}}}. $$ For this
reason, we will focus only on the category ${}^C{\cal YD}(H)_A$ of
right-left quantum Yetter-Drinfel'd modules. It was proved in
Theorem 2.3 and Remark 2.5 of \cite{CMZ1} that ${}^C{\cal YD}(H)_A$
is a special case of the category of Doi-Koppinen modules:
more precisely, if $G=(H,A,C)$ is a Yetter-Drinfel'd datum, then
${\cal O}=(H\ot H^{{\rm op}}, A, C)$ is a Doi-Koppinen datum where $A$
is a left $H\ot H^{{\rm op}}$-comodule algebra via
\begin{equation}\eqlabel{:10}
a\longrightarrow \sum \Bigl( a_{<-1>}\ot S^{-1}(a_{<1>}) \Bigl)
\ot a_{<0>}
\end{equation}
for all $a\in A$ and $C$ is a right $H\ot H^{{\rm op}}$-module
coalgebra via
\begin{equation}\eqlabel{:11}
c\bullet (h\ot k)=k\cdot c\cdot h
\end{equation}
for all $c\in C$, $h$, $k\in H$. Then there exists an isomorphism
of categories $$ {}^C{\cal YD}(H)_A\cong {}^C{\cal M}(H\ot H^{{\rm
op}})_A. $$ Now, we will prove the converse: the Doi-Koppinen
modules category is also a special case of the quantum
Yetter-Drinfel'd category, i.e. both categories are on the same
level of generality. Let ${\cal O}=(H, A, C)$ be a Doi-Koppinen
datum, that is $A$ is a left $H$-comodule algebra and $C$ is a
right $H$-module coalgebra. We view $A$ as an $H$-bicomodule
algebra, where the right $H$-coaction on $A$ is trivial, that is
$A\to A\ot H$, $a\to a\ot 1_H$ for all $a\in A$, and $C$ as an
$H$-bimodule coalgebra where the left action of $H$ on $C$ is also
trivial, that is $H\ot C\to C$, $h\ot c \to \varepsilon (h)c$, for
all $h\in H$, $c\in C$. With these structures, we can view
$G=(H,A,C)$ as a Yetter-Drinfel'd datum and the compatibility
condition \eqref{:2} becomes exactly \eqref{1.1a}, i.e. ${}^C{\cal
M}(H)_A$ is also a special case of ${}^C{\cal YD}(H)_A$.\\ Our
special interest will be corresponding for the case $C=H$. For
this, ${}^H{\cal YD}(H)_A$ will be simply denoted by ${}^H{\cal
YD}_A$, for an arbitrary $H$-bicomodule algebra $A$. An object in
this category is a $k$-module $M$ that is a right $A$-module and a
left $H$-comodule such that
\begin{equation}\eqlabel{:2a}
\rho_M(ma)=\sum S^{-1}(a_{<1>})m_{<-1>}
a_{<-1>}\otimes m_{<0>}a_{<0>}
\end{equation}
for all $m\in M$, $a\in A$. Now, if the right $H$-comodule
structure on $A$ is trivial, that is $A\to A\ot H$, $a\to a\ot 1_H$,
then ${}^H{\cal YD}_A={}^H{\cal M}_A$, the category of relative Hopf
modules.
As ${}^H{\cal M}_A$ is obtained from ${}^H{\cal YD}_A$
by trivializing the right coaction of $H$ on $A$,
we can view the category of quantum Yetter-Drinfel'd modules
${}^H{\cal YD}_A$ as a {\em quantization} of the category of
relative Hopf modules ${}^H{\cal M}_A$.

\section{Total integrals of a Doi-Koppinen datum}\selabel{2}
The point of view expressed in \thref{1.2.12}, evidencing the fact
that integrals in the sense of Doi (or classic integrals on Hopf
algebras) are necessary and sufficient tools for constructing a
natural transformation, leads us to the correct definition of the
integrals for a Doi-Koppinen datum. \\ Let $M\in {}^C{\cal
M}(H)_A$ with the $C$-coaction $\rho_M :M\to C\ot M$. We have seen
in \seref{1.1} that $\rho_M$ is a morphism in ${}^C{\cal M}(H)_A$,
in particular in ${}^C{\cal M}$. We regarded $\rho$ as a natural
transformation between the functors $$ \rho:F_A\circ 1_{{}^C{\cal
M}(H)_A} \to F_A\circ (C \ot \bullet)\circ F^C. $$ Now, in the
light of the above interpretation, a total integral for a
Doi-Koppinen datum should be the necessary and sufficient tool for
constructing a split of $\rho$.

\begin{definition}\delabel{2.1.4a}
Let $(H, A, C)$ be a Doi-Koppinen datum. A $k$-linear map
$\gamma : C\to \Hom(C, A)$ is called an integral of $(H, A, C)$
if:
\begin{eqnarray}\eqlabel{TI1}
\sum c_{(1)}\ot \gamma(c_{(2)})(d)=
\sum d_{(2)}\, \{ \gamma(c)(d_{(1)}) \}_{<-1>}\ot
\{ \gamma(c)(d_{(1)}) \}_{<0>}
\end{eqnarray}
for all $c$, $d\in C$. An integral $\gamma : C\to \Hom(C, A)$
is called total if
\begin{eqnarray}\eqlabel{TI2}
\sum \gamma(c_{(1)})(c_{(2)})=\varepsilon(c)1_A
\end{eqnarray}
for all $c\in C$.
\end{definition}

The concept of integral presented above is obtained by relaxing
the notion of $A$-integral of a Doi-Koppinen datum introduced in
Definition 2.6 of \cite{CIMZ}, leaving aside the $A$-centralizing
condition. More precisely, an $A$-integral is a total integral
$\gamma : C\to \Hom(C, A)$ satisfying the $A$-centralising
condition $$ \sum a_{<0>}\gamma (c a_{<-2>})(da_{<-1>})= \gamma
(c)(d)a $$ for all $a\in A$ and $c$, $d\in C$.

The condition \eqref{TI2} is a normalizing condition: it can be
viewed as a counterpart of the condition $\varphi (1_H)=1_A$,
corresponding to the case $C=H$.

The condition \eqref{TI1} looks very far away from the
colinearity condition that appears in the case $C=H$. However,
if $C$ is projective over $k$, the condition \eqref{TI1}
is in fact a colinearity condition.
Let $c^*\in C^*$. Applying $c^*$ to the first position we obtain
\begin{equation}\eqlabel{2.1.1.4d}
\sum\lan c^*,c_{(1)}\ran\gamma(c_{(2)})(d)=
\sum\lan c^*, d_{(2)}\, \{\gamma(c)(d_{(1)})\}_{<-1>} \ran \,
\{\gamma(c)(d_{(1)})\}_{<0>}
\end{equation}
or equivalent, using \eqref{1.10c}
$$
\gamma (c\lhua c^*)(d)=(\gamma (c)\cdot c^*)(d)
$$
which means that $\gamma$ is right $C^*$-linear.

Furthermore, if $C$ is projective over $k$,
then \eqref{TI1} is equivalent to the fact that $\gamma$ is a
right $C^*$-linear map. In this case, we shall go further:
let $\Hom (C,A)\in {\cal M}_{C^*}$ with the structure from
\eqref{1.10c}. We define
$$
\HOM (C,A)=\Hom (C,A)^{\rm rat}
$$
the rational part of the right $C^*$-module $\Hom(C,A)$.
There is another equivalent way for the definition of
$\HOM (A,C)$: it is the pull-back of the following two maps:
$$
\alpha :\Hom (C,A)\to \Hom (C, C\ot A), \quad
\alpha (f)(c)=\sum
c_{(2)}\cdot f\bigl(c_{(1)}\bigr)_{<-1>} \ot f\bigl(c_{(1)}\bigr)_{<0>}
$$
and
$$
\theta : C\ot \Hom(C, A)\to \Hom (C, C\ot A), \quad
\theta (c\ot f)(d)=c\ot f(d)
$$
for all $c$, $d\in C$ and $f\in \Hom(C,A)$.\\
Being rational as a right $C^*$-module, $\HOM (C, A)$ has a
natural structure of left $C$-comodule (\cite{A}, \cite{S}).
By definition, a $k$-linear map $f: C\to A$ belongs to
$\HOM (C, A)$ if and only if there exists a family of elements
$c_1, \cdots, c_n \in C$ and $h_1, \cdots, h_n \in \Hom(C, A)$
such that
$$
f\cdot c^* = \sum_{i=1}^n h_i \lan c^*, c_i \ran
$$
which is equivalent to
$$
\sum f\bigl(d_{(1)}\bigr)_{<0>}\lan c^*,
d_{(2)}\cdot f\bigl(d_{(1)}\bigr)_{<-1>}\ran
=\sum_{i=1}^n h_i(d) \lan c^*, c_i \ran
$$
for all $c^*\in C^*$, $d\in C$. As $C$ is projective, the last
equation is equivalent to
$$
\sum_{i=1}^n c_i \ot h_i(d)=\sum
d_{(2)}\cdot f\bigl(d_{(1)}\bigr)_{<-1>} \ot f\bigl(d_{(1)}\bigr)_{<0>}
$$
for all $d\in C$.\\
Let now $\gamma : C\to \Hom(C, A)$ be an integral of $(H, A, C)$ and
$f=\gamma(c)$, $c\in C$. Then choosing $c_i=c_{(1)}$ and
$h_i=\gamma (c_{(2)})$ the condition \eqref{TI1} assures that
$\gamma(c) \in \HOM(C,A)$, i.e.
${\rm Im}(\gamma)\subseteq \HOM(C,A)$.

We record these observations in the following

\begin{proposition}
Let $(H, A, C)$ be a Doi-Koppinen datum such that $C$ is projective over
$k$. The following statements are equivalent:
\begin{enumerate}
\item there exists $\gamma : C\to \Hom(C, A)$ an integral of $(H, A, C)$;
\item there exists $\tilde{\gamma}: C\to \HOM(C, A)$ a left $C$-comodule
map.
\end{enumerate}
\end{proposition}

\begin{remarks}
\rm
1. We shall point out now that the integrals introduced above are
in line with Doi's total integrals and with the classic integrals on
Hopf algebras.
For this, let $C=H$ and $\gamma :H\to \Hom(H, A)$ be a total
integral for $(H, A, H)$. For $c=h\in H$ and $d=1_H$, the equation
\eqref{TI1} takes the form
$$
\sum h_{(1)}\ot \gamma (h_{(2)})(1_H)=\sum
\gamma (h)(1_H)_{<-1>} \ot \gamma (h)(1_H)_{<0>}
$$
i.e the map
\begin{equation}\eqlabel{zof1}
\varphi=\varphi_{\gamma}:H\to A, \quad
\varphi(h)=\gamma (h)(1_H)
\end{equation}
for all $h\in H$ is left $H$-colinear, hence a total integral. \\
Conversely, assume that $\varphi :H\to A$ is a total integral and
that the antipode $S$ is bijective (this assumption is given by
the choice of sides (right-left), according to 2) of Remark \ref{zof}). Then
\begin{equation}\eqlabel{zof2}
\gamma:H\to \Hom(H, A), \quad
\gamma (h)(g)=\varphi \Bigl( S^{-1}(g)h\Bigl)
\end{equation}
for all $g$, $h\in H$ is a total integral for $(H, A, H)$. Indeed, for
$c$, $d\in H$ we have
\begin{eqnarray*}
\sum d_{(2)}\gamma(c)(d_{(1)})_{<-1>}\ot \gamma(c)(d_{(1)})_{<0>}
&=& \sum d_{(2)}\varphi \Bigl( S^{-1}(d_{(1)})c\Bigl)_{<-1>}\ot
\varphi \Bigl( S^{-1}(d_{(1)})c\Bigl)_{<0>}\\
(\varphi~{\rm is~left~}H-{\rm colinear})~~
&=& \sum d_{(2)} S^{-1}(d_{(1)})_{(1)} c_{(1)}\ot
\varphi \Bigl( S^{-1}(d_{(1)})_{(2)}c_{(2)}\Bigl) \\
&=& \sum d_{(3)} S^{-1}(d_{(2)}) c_{(1)}\ot
\varphi \Bigl( S^{-1}(d_{(1)}) c_{(2)}\Bigl) \\
&=& \sum c_{(1)} \ot \varphi \Bigl( S^{-1}(d) c_{(2)}\Bigl) \\
&=& \sum c_{(1)} \ot \gamma (c_{(2)})(d)
\end{eqnarray*}
hence $\gamma$ is a total integral in our sense.\\
In particular, at the level of Hopf algebras over a field $k$,
the existence of
an integral of $(H, k, H)$ is equivalent to the existence of
a classical integral on $H$. The correspondence is given by the
formulas \eqref{zof1} and \eqref{zof2}.
We mention that if there exists an integral $\varphi :H\to k$
on $H$, then the antipode $S$ is bijective, hence the formula
\eqref{zof2} can be used.\\
2. The above general definition has an extra-bonus: it
leads to the notion of an integral for
a coalgebra $C$,
\footnote{For $C$ $k$-projective, the notion was introduced
in Definition 4.1 of \cite{CMZ3}.}
which corresponds to the Doi-Koppinen datum
$(k, k, C)$. More precisely, an {\em integral for a coalgebra}
$C$ is a $k$-linear map $\gamma :C\to C^*$ satisfying  the condition
\begin{eqnarray}
\sum c_{(1)}\ot \gamma(c_{(2)})(d)&=&
\sum d_{(2)}\ot \gamma(c)(d_{(1)})
\eqlabel{TI1co}\\
\end{eqnarray}
for all $c$, $d\in C$. Furthermore, $\gamma$ is called
{\em total} if
\begin{eqnarray}
\sum \gamma(c_{(1)})(c_{(2)})&=&\varepsilon(c) \eqlabel{TI2co}
\end{eqnarray}
for all $c$. For a coalgebra $C$ over a field $k$, there exists a
total integral $\gamma :C\to C^*$ if and only if $C$ is coseparable
(Theorem 4.3 of \cite{CMZ3}), and this is the Maschke theorem for
coalgebras.
\end{remarks}

We have proved that classical integrals or Doi's total
integrals are examples of integrals in our sense.
We shall indicate now two important classes of examples that are
not produced in this way.

\begin{examples}
\rm
1. Let $C={\rm M}^n(k)$ be the
$n\times n$ comatrix coalgebra, i.e. $C$ is the coalgebra
with a $k$-basis $\{c_{ij}\mid i,j=1,\cdots,n \}$ such that
$$
\Delta(c_{ij})=\sum_{u=1}^{n}c_{iu}\ot c_{uj}, \quad
\varepsilon (c_{ij})=\delta_{ij}
$$
for all $i$, $j=1,\cdots, n$.
Let $(H, A, C)=(H, A, {\rm M}^n(k))$ be a Doi-Koppinen datum and
$B=A^{{\rm co}(H)}=\{\; a\in A \;\mid \; \rho(a)=1_H\ot a \;\}$
be the subalgebra of coinvariants of $A$. Let
$\mu=(\mu_{ij}) \in {\rm M}_n(B)$ be an arbitrary
$n\times n$-matrix over $B$ (for example,
$\mu_{ij}=a_{ij}1_A$, where $a_{ij}$ are scalars of $k$).
Then the map
$$
\gamma=\gamma_{\mu}:\ C\to \Hom(C, A), \quad
\gamma (c_{ij})(c_{rs})= \delta_{is}\mu_{rj}
$$
is an integral of $(H, A, {\rm M}^n(k))$. Indeed, for
$c=c_{ij}$ and $d=c_{kl}$
$$
\sum c_{(1)}\ot \gamma(c_{(2)})(d)=
\sum_{u=1}^{n} c_{iu}\ot \delta_{ul}\mu_{kj}=
c_{il}\ot \mu_{kj}
$$
and
\begin{eqnarray*}
\sum d_{(2)}\gamma(c)(d_{(1)})_{<-1>}\ot \gamma(c)(d_{(1)})_{<0>}
&=&\sum_{v=1}^{n} c_{vl} \gamma (c_{ij})(c_{kv})_{<-1>}\ot
\gamma (c_{ij})(c_{kv})_{<0>}\\
(~\mu_{kj}\in B~)
&=&\sum_{v=1}^{n} c_{vl}\ot \delta_{iv}\mu_{kj}=
c_{il}\ot \mu_{kj}
\end{eqnarray*}
i.e. $\gamma_{\mu}$ is an integral. On the other hand,
for $c=c_{ij}$,
$$
\sum \gamma (c_{(1)})(c_{(2)})=
\sum_{u=1}^{n} \gamma (c_{iu})(c_{uj})=
\sum_{u=1}^{n}\delta_{ij} \mu_{uu}=
\delta_{ij}{\rm Tr}(\mu)
$$
Hence, $\gamma_{\mu}$ is a total integral if and only if Tr$(\mu)=1_A$. \\
2. Another class of examples of total integrals arises from the
graded case. Let $X$ be a set and $C=kX$ be the group-like
coalgebra, i.e. $kX$ is the free $k$-module having $X$ as a basis
and $\Delta (x)=x\ot x$, $\varepsilon (x)=1$, for all $x\in X$.
Let $(H, A, C)=(H, A, kX)$ be a Doi-Koppinen datum and
$\mu=(\mu_{xy})_{x, y\in X}$ be a family of elements of
$B=A^{{\rm co}(H)}$. The map
$$
\gamma=\gamma_{\mu} :kX\to \Hom(kX, A), \quad
\gamma (x)(y)=\delta_{xy}\mu_{xy}
$$
for all $x$, $y\in X$ is an integral of $(H, A, kX)$. Indeed,
for $x$, $y\in X$ we have
$$
\sum y\gamma(x)(y)_{<-1>}\ot \gamma(x)(y)_{<0>}=
y\ot \delta_{xy}\mu_{xy} =x\ot \delta_{xy}\mu_{xy}=
x\ot \gamma(x)(y)
$$
i.e. $\gamma_{\mu}$ is an integral of $(H, A, kX)$. Furthermore,
$\gamma_{\mu}$ is a total integral iff $\mu_{xx}=1_A$ for all $x\in X$. \\
The category ${}^{kX}{\cal M}(H)_A$ associated to this Doi-Koppinen datum
covers a large class of examples of $X$-graded representations of $A$,
starting from the category of super-graded vector spaces
(corresponding to the trivial case $H=A=k$, $k$ a field) to the
category of graded modules by $G$-sets (corresponding to the case
$H=kG$, where $G$ is a group acting on the set $X$).
\end{examples}

We shall now prove that the existence of an integral
$\gamma :C\to \Hom(C, A)$ permits
the deformation of a $k$-linear map between two Doi-Koppinen modules
until it becomes a $C$-colinear map.

\begin{proposition} \prlabel{greva}
Let $(H,A,C)$ be a Doi-Koppinen datum, $M\in {}^C{\cal M}(H)_A$,
$N\in {}^C{\cal M}$ and $u:N\to M$ a $k$-linear map.
Suppose that there exists
$\gamma : C\to \Hom(C, A)$ an integral. Then:
\begin{enumerate}
\item the map
$$
\tilde{u}: N\to M, \quad
\tilde{u} (n)= \sum
u(n_{<0>})_{<0>}\gamma (n_{<-1>})(u(n_{<0>})_{<-1>})
$$
for all $n\in N$, is left $C$-colinear;
\item if $\gamma$ is a total integral and
$f: M\to N$ is a morphism in ${}^C{\cal M}(H)_A$ which is a
$k$-split injection (resp. a $k$-split surjection), then $f$ has
a $C$-colinear retraction (resp. a section).
\end{enumerate}
\end{proposition}

\begin{proof}
1. For $n\in N$ we have
\begin{eqnarray*}
\rho_M (\tilde{u}(n)) &=&
\sum u(n_{<0>})_{<-1>} \gamma (n_{<-1>})
\Bigl (u(n_{<0>})_{<-2>}\Bigl)_{<-1>}\\
&&\ot u(n_{<0>})_{<0>} \gamma (n_{<-1>})
\Bigl (u(n_{<0>})_{<-2>}\Bigl)_{<0>} \\
&=&\sum u(n_{<0>})_{<-1>(2)} \gamma (n_{<-1>})
\Bigl (u(n_{<0>})_{<-1>(1)}\Bigl)_{<-1>}\\
&&\ot u(n_{<0>})_{<0>} \gamma (n_{<-1>})
\Bigl (u(n_{<0>})_{<-1>(1)}\Bigl)_{<0>} \\
\eqref{TI1}~~~&=&
\sum n_{<-2>} \ot u(n_{<0>})_{<0>}
\gamma (n_{<-1>})(u(n_{<0>})_{<-1>})\\
&=&\sum n_{<-1>} \ot \tilde{u}(n_{<0>})\\
&=& ({\rm Id}\ot \tilde{u})\rho_{N}(n)
\end{eqnarray*}
hence $\tilde{u}$ is left $C$-colinear.\\
2. Let $u:N\to M$ be a $k$-linear retraction (resp. section) of $f$.
Then $\tilde{u}: N\to M$ is a left $C$-colinear retraction
(resp. section) of $f$. Assume first that $u$ is a retraction of $f$.
Then, for $m\in M$
\begin{eqnarray*}
(\tilde{u}\circ f)(m) &=&\sum
u(f(m)_{<0>})_{<0>}\gamma (f(m)_{<-1>})(u(f(m)_{<0>})_{<-1>})\\
(~{\rm f~is~}C-{\rm colinear})~~
&=&\sum
u(f(m_{<0>}))_{<0>}\gamma (m_{<-1>})(u(f(m_{<0>}))_{<-1>})\\
(~u\circ f={\rm Id}~)
&=& \sum m_{<0>}\gamma(m_{<-2>})(m_{<-1>})=m
\end{eqnarray*}
hence $\tilde{u}: N\to M$ is a left $C$-colinear retraction of $f$.
On the other hand, if $u$ is a section of $f$, then for $n\in N$
\begin{eqnarray*}
(f\circ \tilde{u})(n) &=&\sum
f\Bigl (u(n_{<0>})_{<0>}\gamma (n_{<-1>})(u(n_{<0>})_{<-1>})\Bigl)\\
(~{\rm f~is~}A-{\rm linear})~~
&=&\sum
f\Bigl (u(n_{<0>})_{<0>} \Bigl)
\gamma (n_{<-1>})(u(n_{<0>})_{<-1>})\\
(~{\rm f~is~}C-{\rm colinear})~~
&=&\sum
\bigl (f(u(n_{<0>}))\bigl )_{<0>}\gamma (n_{<-1>})
\Bigl(\bigl (f(u(n_{<0>}))\bigl )_{<-1>}\Bigl)\\
(~f\circ u={\rm Id}~)
&=& \sum n_{<0>}\gamma(n_{<-2>})(n_{<-1>})=n
\end{eqnarray*}
i.e. $\tilde{u}: N\to M$ is a left $C$-colinear section of $f$.
\end{proof}

We will prove now the version of \thref{1.2.12} for Doi-Koppinen modules.
We note that the counterpart of the item 3) of \thref{1.2.12}
has a different form: the difference is given by the fact that
for $C=H$, $A\in {}^H{\cal M}_A$ with the natural structures,
while for an arbitrary $C$, $A$ does not have a structure
of object in ${}^C{\cal M}(H)_A$. Thus, $A$ is replaced now with
$C\ot A$ and this time, $\rho_{C\ot A}^l$ splits $C$-bicolinearly.
Parts of the proof of the following theorem are closely
related to the ideas presented in Section 2 of \cite{CIMZ}.

\begin{theorem}\thlabel{doi}
Let $(H,A,C)$ be a Doi-Koppinen datum. The following statements are
equivalent:
\begin{enumerate}
\item there exists $\gamma : C\to \Hom(C, A)$ a total integral;
\item the natural transformation
$\rho:F_A\circ 1_{{}^C{\cal M}(H)_A} \to
F_A\circ (C \ot \bullet)\circ F^C$ splits;
\item the left $C$-coaction on $C\ot A$,
$\;\rho_{C\ot A}^l :C\ot A\to C\ot C\ot A, \;
\rho_{C\ot A}^l (c\ot a)= \sum c_{(1)}\ot c_{(2)}\ot a$
splits in ${}^C{\cal M}^C$, the category of $C$-bicomodules.
\end{enumerate}
Consequently, if one of the equivalent conditions holds,
any Doi-Koppinen module is relative injective as a left $C$-comodule.
\end{theorem}

\begin{proof}
$1) \Rightarrow 2) $ Let $\gamma: C\to \Hom (C, A)$ be a total integral.
We have to construct a natural transformation $\lambda$ that
splits $\rho$. Let $M \in {}^C{\cal M}(H)_A$ and
$u_M: C\ot M\to M$, be the $k$-linear retraction of
$\rho_M:M \to C\ot M$ given by
$u_M (c\ot m)= \varepsilon(c)m$, for all $c\in C$ and $m\in M$.
We define $\lambda_M=\tilde{u_M}$, i.e.
\begin{eqnarray}\eqlabel{mm13}
\lambda_M=\lambda_M(\gamma): C\ot M\to M, \quad
\lambda_M (c\ot m)=\sum m_{<0>} \gamma (c)(m_{<-1>})
\end{eqnarray}
for all $c$, $m\in M$.
It follows from \prref{greva} that the map
$\lambda_M$ is a left $C$-colinear retraction of $\rho_M$.

It remains to prove that $\lambda=\lambda(\gamma)$ is a natural
transformation. Let $f: M\to N$ be a morphism in
${}^C{\cal M}(H)_A$. We have to prove that the diagram
$$
\begin{diagram}
C\ot M&\rTo{\lambda_M}&M\\
\dTo{{\rm Id}\ot f}&&\dTo_{f}\\
C\ot N&\rTo{\lambda_N}&N\\
\end{diagram}
$$
is commutative. For $c\in C$ and $m\in M$, using that $f$ is
right $A$-linear, we have
$$
(f\circ \lambda_M) (c\ot m)=\sum f(m_{<0>} \gamma (c)(m_{<-1>}))=
\sum f(m_{<0>}) \gamma (c)(m_{<-1>})
$$
and using that $f$ is left $C$-colinear
\begin{eqnarray*}
\Bigl (\lambda_N \circ ({\rm Id}\ot f) \Bigl) (c\ot m)&=&
\lambda_N(c\ot f(m))\\
&=& \sum  f(m)_{<0>} \gamma (c)(f(m)_{<-1>})\\
&=& \sum f(m_{<0>}) \gamma (c)(m_{<-1>})
\end{eqnarray*}
i.e. $\lambda$ is a natural transformation that splits $\rho$.\\
$2) \Rightarrow 3) $ Assume that for any $M\in {}^C{\cal M}(H)_A$
the $C$-coaction $\rho_M :M\to C\ot M$ splits in the category
${}^C{\cal M}$ of left $C$-comodules and the character of the
splitting is functorial. In particular,
$\rho_{C\ot A}^l :C\ot A\to C\ot C\ot A$ splits in ${}^C{\cal M}$,
and let $\lambda=\lambda_{C\ot A}: C\ot C\ot A\to C\ot A$ be a
left $C$-colinear retraction of it. Using the naturality of
$\lambda$,
we will prove that $\lambda$ is also right $C$-colinear, where
$C\ot A$ and $C\ot C\ot A$ are right $C$-comodules via
\eqref{sase} and \eqref{11}.\\
First, let $V$ be a $k$-module and $M\in {}^C{\cal M}(H)_A$. Then
$M\ot V\in {}^C{\cal M}(H)_A$ via the structures arising
from the ones of $M$, i.e
$$
(m\ot v)a= ma\ot v, \quad
\rho_{M\ot V}=\rho_M \ot {\rm Id}_V
$$
for all $m\in M$, $a\in A$ and $v\in V$. Using the naturality
of $\lambda$, we shall prove now that
\begin{equation}\eqlabel{mm15}
\lambda_{M\ot V}=\lambda_{M}\ot {\rm Id}_V
\end{equation}
Let $v\in V$ and
$$
g_v :M\to M\ot V, \quad
g_v(m)= m \ot v
$$
Then $g_v$ is a morphism in ${}^C{\cal M}(H)_A$. From the
naturality of $\lambda$ we obtain that
$$
g_v \circ \lambda_{M}=\lambda_{M\ot V} \circ ({\rm Id}_C \ot g_v)
$$
Hence
$$
\lambda_{M\ot V} (c\ot m\ot v)=
g_v(\lambda_{M}(c\ot m))=
\lambda_{M}(c\ot m)\ot v
$$
i.e. \eqref{mm15} holds. In particular, let us take
$M=C\ot A$ and $V=C$ viewed only as a $k$-module. Then
$C\ot A\ot C\in {}^C{\cal M}(H)_A$ via the structures arising
from the ones of $C\ot A$, i.e
\begin{eqnarray}
(c\ot a\ot d)b&=& \sum cb_{<-1>}\ot ab_{<0>}\ot d\eqlabel{mm14}\\
\rho_{C\ot A\ot C}(c\ot a\ot d)&=&
\sum c_{(1)}\ot c_{(2)}\ot a\ot d\eqlabel{mm141}
\end{eqnarray}
for all $c$, $d\in C$, $a$, $b\in A$. With these structures the map
$$
f=\rho^{r}_{C\ot A}: C\ot A \to C\ot A\ot C, \quad
f(c\ot a)=\sum c_{(1)}\ot a_{<0>}\ot c_{(2)}S(a_{<-1>})
$$
is a morphism in ${}^C{\cal M}(H)_A$. From the
naturality of $\lambda$ the following diagram
\begin{equation}\eqlabel{mm16}
\begin{diagram}
C\ot C\ot A&\rTo{\lambda_{C\ot A}}&C\ot A\\
\dTo{I_C\ot f=\rho_{C\ot C\ot A}^{r}}&&\dTo_{f=\rho_{C\ot A}^{r}}\\
C\ot C\ot A\ot C&\rTo{\lambda_{C\ot A\ot C}=
\lambda_{C\ot A}\ot I_C}&C\ot A\ot C
\end{diagram}
\end{equation}
is commutative, i.e. $\lambda=\lambda_{C\ot A}$ is also right
$C$-colinear.\\
$3) \Rightarrow 1) $ The left $C$-coaction
$\rho^{l}_{C\ot A}: C\ot A\to C\ot C\ot A$ is a $C$-bicomodule map.
Let $\lambda=\lambda_{C\ot A}: C\ot C\ot A\to C\ot A$ be a
split of $\rho^{l}_{C\ot A}$ in ${}^C{\cal M}^C$. In particular,
$$
\lambda (\sum c_{(1)}\ot c_{(2)}\ot a)=c\ot a
$$
for all $c\in C$, $a\in A$. We define,
\begin{equation}\eqlabel{mmti}
\gamma : C\to \Hom (C, A), \quad
\gamma (c)(d):= (\varepsilon \ot {\rm Id}) \lambda (c\ot d\ot 1_A)
\end{equation}
for all $c$, $d\in C$. We will prove that $\gamma$ is a total integral.
First,
$$
\sum \gamma (c_{(1)})(c_{(2)})= \sum
(\varepsilon \ot {\rm Id}_A) \lambda (c_{(1)}\ot c_{(2)}\ot 1_A)=
\varepsilon (c)1_A
$$
i.e. \eqref{TI2} holds. We will prove that \eqref{TI1} also holds.
For $c$, $d\in C$, the left hand side of \eqref{TI1} is
\begin{eqnarray*}
\sum c_{(1)}\ot \gamma(c_{(2)}) (d)&=& \sum c_{(1)}\ot
(\varepsilon \ot {\rm I}_A)\lambda (c_{(2)}\ot d\ot 1_A)\\
(\lambda~{\rm is~left~}C-{\rm colinear})~~ &=&({\rm Id}_C \ot
\varepsilon \ot {\rm Id}_A) \rho_{C\ot A}\Bigl (\lambda (c\ot d\ot
1_A)\Bigl)\\ &=&\lambda (c\ot d\ot 1_A)
\end{eqnarray*}
In order to compute the right hand side of \eqref{TI1} we adopt the
temporary notation
$$\lambda (c\ot d_{(1)}\ot 1_A)=\sum_i c_i\ot a_i.$$
Now,
\begin{eqnarray*}
\sum d_{(2)}\gamma(c)(d_{(1)})_{<-1>}\ot \gamma(c)(d_{(1)})_{<0>}&=&
\sum d_{(2)}\varepsilon (c_{i})a_{i_{<-1>}}\ot a_{i_{<0>}}\\
&=& \sum d_{(2)}a_{i_{<-1>}}\ot \varepsilon (c_{i})a_{i_{<0>}}
\end{eqnarray*}
Hence, the equation \eqref{TI1} is equivalent to
\begin{eqnarray}\eqlabel{2.1.1.1b}
\lambda(c\ot d\ot 1_A)=
\sum (d_{(2)}\ot 1_A)
(\varepsilon_C\ot \rho_A)\lambda (c\ot d_{(1)}\ot 1_A)
\end{eqnarray}
for all  $c$, $d\in C$. Now, it is time to use the fact
that $\lambda$ is also right $C$-colinear. Denoting
$$
\lambda (c\ot d\ot 1_A)=\sum_i D_i\ot A_i \in C\ot A
$$
and evaluating the diagram \eqref{mm16} at
$c\ot d \ot 1_A$, we obtain
$$\sum \lambda (c\ot d_{(1)}\ot 1_A)\ot d_{(2)}=
\sum D_{i_{(1)}}\ot A_{i_{<0>}}\ot D_{i_{(2)}}S(A_{i_{<-1>}})$$
hence,
$$\sum d_{(2)}\ot (\varepsilon_C\ot I_A)\bigl(\lambda
(c\ot d_{(1)}\ot 1_A)\bigr)=
\sum D_{i}S\bigl(A_{i_{<-1>}}\bigr)\ot A_{i_{<0>}}$$
Now we apply $\rho_A$ to the second factor of both sides. Using the
fact that $\rho_A\circ (\varepsilon_C\ot I_A)=\varepsilon_C\ot \rho_A$,
we obtain
$$\sum d_{(2)}\ot (\varepsilon_C\ot \rho_A)
\bigl(\lambda(c\ot d_{(1)}\ot 1_A)\bigr)=
\sum D_{i}S\bigl(A_{i_{<-2>}}\bigr)\ot A_{i_{<-1>}}\ot A_{i_{<0>}}.$$
\eqref{2.1.1.1b} follows after we let the second factor act on the
first one.
\end{proof}

Leaving aside the normalizing condition \eqref{TI2}, we obtain the
following directly from the proof of the theorem:

\begin{corollary}
Let $(H,A,C)$ be a Doi-Koppinen datum. The following statements are
equivalent:
\begin{enumerate}
\item there exists $\gamma : C\to \Hom(C, A)$ an integral of
$(H, A, C)$;
\item there exists
$\lambda: F_A\circ (C \ot \bullet)\circ F^C \to
F_A\circ 1_{{}^C{\cal M}(H)_A}$
a natural transformation;
\item there exists $\lambda^{\prime}: C\ot C\ot A \to C\ot A$
a $C$-bicomodule map .
\end{enumerate}
\end{corollary}

\begin{remark} \rm
The \thref{doi} has an interesting consequence in the finite
dimensional case.
First, let $R\subset S$ be a ring extension. $S/R$ is called a
RIT-{\em extension} (right integral type) if any right $S$-module
is injective as a right $R$-module. Assume that $C$ is finite
dimensional over a field $k$ and $\gamma : C\to \Hom (C, A)$ is
a total integral. Then $C^* \subset A\# C^*$ is a RIT-extension.
Indeed, as $C$ is finite dimensional, the functor
$P:\ {}^C{\cal M}(H)_A\to {\cal M}_{A\# C^*}$ is an equivalence
of categories, which allows us to apply \thref{doi}.
\end{remark}

The object $C\ot A$ plays a special role in ${}^C{\cal M}(H)_A$:
first, if $k$ is a field, it was proved in \cite{CR} that $C\ot A$
is a subgenerator in ${}^C{\cal M}(H)_A$, i.e. any object in
${}^C{\cal M}(H)_A$ is isomorphic to a subobject of a quotient of
direct sums of copies of $C\ot A$. Secondly, over a commutative ring,
Corollary 2.9 of \cite{CMZ2} proves that, if the
forgetful functor $F^C: {}^C{\cal M}(H)_A\to {\cal M}_A$ is
Frobenius (i.e. by definition has the same left and right adjoint),
then $C\ot A$ is a generator in ${}^C{\cal M}(H)_A$.
Finally, Lemma 2.8 of \cite{DNT} shows that, if $k$ is a field
and $C$ is a left and right quasi co-Frobenius coalgebra,
then $C\ot A$ is a generator in ${}^C{\cal M}(H)_A$.
We shall prove now the main applications of the existence of a total
integral.

\begin{theorem}\thlabel{mm3}
Let $(H,A,C)$ be a Doi-Koppinen datum and suppose that
there exists $\gamma : C\to \Hom(C, A)$ a total integral. Then
for any $M\in {}^C{\cal M}(H)_A$ the map
\begin{equation}\eqlabel{mm18}
f: C\ot A\ot M \to M, \qquad
f(c\ot a\ot m)=\sum m_{<0>} \gamma (c S(a_{<-1>}))(m_{<-1>})a_{<0>}
\end{equation}
for all $c\in C$, $a\in A$ and $m\in M$ is a $k$-split epimorphism
in ${}^C{\cal M}(H)_A$. \\
In particular, $C\ot A$ is a generator in the category
${}^C{\cal M}(H)_A$.
\end{theorem}

\begin{proof}
$C\ot A\ot M$ is viewed as an object in ${}^C{\cal M}(H)_A$ with the
structures arising from the ones of $C\ot A$, i.e
\begin{eqnarray*}
(c\ot a\ot m)b&=& \sum cb_{<-1>}\ot ab_{<0>}\ot m\\
\rho_{C\ot A\ot M}(c\ot a\ot m)&=&
\sum c_{(1)}\ot c_{(2)}\ot a\ot m
\end{eqnarray*}
for all $c\in C$, $a$, $b\in A$ and $m\in M$. First we shall
prove that $f$ is a $k$-split (even a $C$-colinear split) surjection.
Let $g: M\to C\ot A\ot M$, $g(m)=\sum m_{<-1>}\ot 1_A\ot m_{<0>}$,
for all $m\in M$. Then $g$ is left $C$-colinear
(but is not right $A$-linear) and for $m\in M$ we have
\begin{eqnarray*}
(f\circ g)(m)&=&\sum f(m_{<-1>}\ot 1_A\ot m_{<0>})\\
&=&\sum m_{<0><0>}\gamma (m_{<-1>})(m_{<0><-1>})\\
&=&\sum m_{<0>}\gamma (m_{<-2>})(m_{<-1>})\\
&=&\sum m_{<0>}\gamma (m_{<-1>(1)})(m_{<-1>(2)})\\
&=&\sum m_{<0>}\varepsilon (m_{<-1>})=m
\end{eqnarray*}
i.e. $g$ is a left $C$-colinear section of $f$. For $a$, $b\in A$,
$c\in C$ $m\in M$ we have
\begin{eqnarray*}
f((c\ot a\ot m)b)&=& \sum f(cb_{<-1>}\ot ab_{<0>}\ot m)\\
&=&\sum m_{<0>}\gamma \Bigl( cb_{<-1>}S(b_{<0><-1>})S(a_{<-1>})\Bigl)
(m_{<-1>})a_{<0>}b_{<0><0>}\\
&=&\sum m_{<0>}\gamma \Bigl( cb_{<-2>}S(b_{<-1>})S(a_{<-1>})\Bigl)
(m_{<-1>})a_{<0>}b_{<0>}\\
&=&\sum m_{<0>}\gamma \Bigl( cS(a_{<-1>})\Bigl)
(m_{<-1>})a_{<0>}b\\
&=& f(c\ot a\ot m)b
\end{eqnarray*}
i.e. $f$ is right $A$-linear. It remains to prove that $f$ is also
left $C$-colinear. First we compute
\begin{eqnarray*}
(I_C \ot f)\rho_{C\ot A\ot M} (c\ot a\ot m)&=&
\sum c_{(1)}\ot f(c_{(2)}\ot a \ot m)\\
&=& \sum c_{(1)}\ot m_{<0>} \gamma (c_{(2)} S(a_{<-1>}))(m_{<-1>})a_{<0>}
\end{eqnarray*}
and
\begin{eqnarray*}
\rho_M (f(c\ot a\ot m)) &=& \sum
\rho_M (m_{<0>} \gamma (c S(a_{<-1>}))(m_{<-1>})a_{<0>} )\\
&=& \sum m_{<0><-1>} \Bigl( \gamma (c S(a_{<-1>}))(m_{<-1>})\Bigl )_{<-1>}
a_{<0><-1>} \ot \\
&&m_{<0><0>} \Bigl(\gamma (c S(a_{<-1>}))(m_{<-1>})\Bigl )_{<0>}
a_{<0><0>}\\
&=& \sum m_{<-1>} \Bigl( \gamma (c S(a_{<-2>}))(m_{<-2>})\Bigl )_{<-1>}
a_{<-1>} \ot \\
&&m_{<0>} \Bigl(\gamma (c S(a_{<-2>}))(m_{<-2>})\Bigl )_{<0>}
a_{<0>}\\
&=& \sum \underline{
m_{<-1>(2)} \Bigl( \gamma
(c S(a_{<-2>}))(m_{<-1>(1)})\Bigl )_{<-1>} }
a_{<-1>} \ot \\
&&m_{<0>}\underline{
\Bigl(\gamma (c S(a_{<-2>}))(m_{<-1>(1)})\Bigl )_{<0>} }
a_{<0>}\\
\eqref{TI1}~~&=& \sum c_{(1)}S(a_{<-2>})_{(1)}a_{<-1>}\ot \\
&& m_{<0>}\gamma (c_{(2)} S(a_{<-2>})_{(2)})(m_{<-1>})a_{<0>}\\
&=&\sum c_{(1)}S(a_{<-2>})a_{<-1>}\ot \\
&& m_{<0>}\gamma (c_{(2)} S(a_{<-3>}))(m_{<-1>})a_{<0>}\\
&=& \sum c_{(1)} \ot m_{<0>}\gamma (c_{(2)} S(a_{<-1>}))(m_{<-1>})a_{<0>}\\
&=& (I_C\ot f)\rho_{C\ot A\ot M} (c\ot a\ot m)
\end{eqnarray*}
i.e. $f$ is left $C$-colinear. Hence, we proved that $f$ is an
epimorphism in ${}^C{\cal M}(H)_A$ and has a $C$-colinear section.\\
Now, taking a $k$-free presentation of $M$ in the category of
$k$-modules
$$
k^{(I)}\stackrel{\pi}{\longrightarrow}M{\longrightarrow}0
$$
we obtain an epimorphism in ${}^C{\cal M}(H)_A$
$$
(C\ot A)^{(I)}\cong C\ot A\ot k^{(I)}\stackrel{g}
{\longrightarrow}M{\longrightarrow}0
$$
where $g= f\circ (I_C\ot I_A\ot \pi)$. Hence $C\ot A$ is a
generator in ${}^C{\cal M}(H)_A$.
\end{proof}

\section{The affineness criterion for quantum Yetter-Drinfel'd modules}
\selabel{3}
Let $H$ be a Hopf algebra with a bijective antipode and
$A$ an $H$-bicomodule algebra.
Let ${}^H{\cal YD}_A$ be the category of (right-left)
quantum Yetter-Drinfel'd modules, i.e. an object in
${}^H{\cal YD}_A$ is a right $A$-module, left $H$-comodule
$(M, \cdot, \rho_M)$ such that
\begin{equation}\eqlabel{gal11}
\sum m_{<-1>}a_{<-1>}\ot m_{<0>}\cdot a_{<0>}=
\sum a_{<1>} (m\cdot a_{<0>})_{<-1>}\ot (m\cdot a_{<0>})_{<0>}
\end{equation}
for all $m\in M$ and $a\in A$.

\begin{remark}
\rm Let $M$ be a right $A$-module and a left $H$-comodule. Then
the compatibility relation \eqref{gal11} is equivalent to
\begin{equation}\eqlabel{gal1}
\rho_M(m\cdot a)=\sum S^{-1}(a_{<1>})m_{<-1>}a_{<-1>}\ot
m_{<0>}\cdot a_{<0>}
\end{equation}
for all $m\in M$ and $a\in A$. Indeed, assume first
that \eqref{gal1} holds. Then for $a\in A$, $m\in M$
\begin{eqnarray*}
\sum a_{<1>} (m\cdot a_{<0>})_{<-1>}\ot (m\cdot a_{<0>})_{<0>}
&=&\sum a_{<1>} S^{-1}(a_{<0><1>})m_{<-1>}a_{<0><-1>}\\
&&\ot  m_{<0>}\cdot a_{<0><0>}\\
&=&\sum m_{<-1>}a_{<-1>}\ot m_{<0>}\cdot a_{<0>}
\end{eqnarray*}
Conversely, if \eqref{gal11} holds then
\begin{eqnarray*}
\rho_M (m\cdot a)&=&\sum (m\cdot a)_{<-1>}\ot (m\cdot a)_{<0>} \\
&=&\sum \varepsilon (a_{<1>})
(m\cdot a_{<0>})_{<-1>}\ot (m\cdot a_{<0>})_{<0>}\\
&=& \sum  S^{-1}(a_{<1>}) a_{<0><1>}
(m\cdot a_{<0><0>})_{<-1>}\ot (m\cdot a_{<0><0>})_{<0>}\\
&=&\sum S^{-1}(a_{<1>})m_{<-1>}a_{<-1>}\ot
m_{<0>}\cdot a_{<0>}
\end{eqnarray*}
fro all $a\in A$, $m\in M$.
\end{remark}

\begin{examples}
\rm  1. Let $A=H$ and $\rho^l=\rho^r=\Delta$. Then ${}^H{\cal
YD}_H$ is the category of crossed (or Yetter-Drinfel'd) modules
introduced in \cite{Y} (see also \cite{RT} for all left-right
conventions). \\ 2. If $\rho_r :A\to A\ot H$ is the trivial
coaction, that is $\rho^r(a)=a\ot 1_H$, then ${}^H{\cal
YD}_A={}^H{\cal M}_A$, the category of classical relative Hopf
modules.\\ 3. If both coactions $\rho^l$ and $\rho^r$ are trivial,
then ${}^H{\cal YD}_A={}^H{\cal L}_A$, the category of Long
dimodules. This category was defined by F.W. Long in \cite{L} for
the case $A=H$, a commutative and cocommutative Hopf algebra and
was studied in connection with the Brauer group. In the general
case, the category ${}^C{\cal L}_A$ was introduced in \cite{Mi2}
and was studied related to a nonlinear equation. \\ 4. Let $H=kG$
where $G$ is a group. Then a $kG$-bicomodule algebra is a
bi-graded $k$-algebra A, having two compatible gradations of type
$G$ and an object in  ${}^{kG}{\cal YD}_A$ is a $G$-graded
representation on $A$ such that the $A$-action agrees with the
bi-graduation.
\end{examples}

An important object of ${}^H{\cal YD}_A$
is the Verma structure $(A, \cdot, \tilde{\rho})$, where
$\cdot$ is the multiplication on $A$ and the
left $H$-coaction $\tilde{\rho}$ given by
\begin{equation}\eqlabel{gal19}
\tilde{\rho}: A\to H\ot A, \quad
\tilde{\rho} (a)=\sum S^{-1}(a_{<1>})a_{<-1>}\ot a_{<0>}
\end{equation}
for all $a\in A$. Then
$(A, \cdot, \tilde{\rho})\in {}^H{\cal YD}_A$ and we will see that
it will play a crucial role in this section.
In the particular case $(A=H, \rho^l=\rho^r=\Delta)$, the above
structure is just the right-left version of the Verma
Yetter-Drinfel'd module over $H$ defined in equation (2.6) of
\cite{EK}.

\begin{remark}\relabel{galred}
\rm  If $H$ is commutative then $(A, \tilde{\rho})$ is a structure
of left $H$-comodule algebra on $A$. Hence, we can associate the
usual category ${}^H{\cal M}_A$ of classical Hopf modules to it:
it is easy to show that ${}^H{\cal YD}_A={}^H{\cal M}_A$, i.e.
in the commutative case the theory of the category of quantum
Yetter-Drinfel'd modules can be reduced to the study of the Hopf
modules category. The theory presented below is relevant for the
noncommutative case.
\end{remark}

We view now  ${}^H{\cal YD}_A$ as the category of Doi-Koppinen
modules associated to the Doi-Koppinen datum $(H\ot H^{{\rm op}}, A, H)$
where

$\bullet$ $A$ is a left $H\ot H^{{\rm op}}$-comodule algebra via
\begin{equation}\eqlabel{gal10}
a\longrightarrow \sum \Bigl( a_{<-1>}\ot S^{-1}(a_{<1>}) \Bigl)
\ot a_{<0>}
\end{equation}
for all $a\in A$ and

$\bullet$ $H$ is a right $H\ot H^{{\rm op}}$-module coalgebra via
\begin{equation}\eqlabel{gal111}
g\cdot (h\ot k)=kgh
\end{equation}
for all $g$, $h$, $k\in H$. Then
${}^H{\cal M}(H\ot H^{{\rm op}})_A={}^H{\cal YD}_A$, and hence all the
concepts, structures and results from previous sections can
be formulated for ${}^H{\cal YD}_A$.
For instance, $H\ot A\in {}^H{\cal YD}_A$ via the
following structures
\begin{eqnarray}
(h\ot b)a&=& \sum S^{-1}(a_{<1>})h a_{<-1>} \ot ba_{<0>}\eqlabel{gal9}\\
\rho_{H\ot A}(h\ot b)&=& \sum h_{(1)}\ot h_{(2)}\ot b\eqlabel{gal101}
\end{eqnarray}
for all $h\in H$, $a$, $b\in A$. The \deref{2.1.4a} has
the following form

\begin{definition}\delabel{galint}
Let $H$ be a Hopf algebra with a bijective antipode and $A$ an
$H$-bicomodule algebra. A $k$-linear map
$\gamma : H\to \Hom(H, A)$ is called a quantum integral
if:
\begin{eqnarray}\eqlabel{galTI1}
\sum g_{(1)}\ot \gamma(g_{(2)})(h)=
\sum S^{-1} \Bigl (\{\gamma (g)(h_{(1)})\}_{<1>}\Bigl)
h_{(2)}\{\gamma(g)(h_{(1)})\}_{<-1>}\ot \{\gamma(g)(h_{(1)})\}_{<0>}
\end{eqnarray}
for all $g$, $h\in H$. A quantum integral $\gamma: H\to \Hom(H, A)$
is called total if
\begin{eqnarray}\eqlabel{galTI2}
\sum \gamma(h_{(1)})(h_{(2)})=\varepsilon(h)1_A
\end{eqnarray}
for all $h\in H$.
\end{definition}

\begin{remarks}\relabel{galoi}
\rm 1. Let $\gamma : H\to \Hom(H, A)$ be a quantum integral.
Then,
$$
\varphi=\varphi_{\gamma}: H\to A, \quad
\varphi (g)=\gamma (g)(1_H)
$$
satisfies the condition
$$
\sum g_{(1)} \ot \varphi(g_{(2)}) =
\sum S^{-1} \Bigl (\varphi (g)_{<1>}\Bigl)
\varphi(g)_{<-1>}\ot \varphi(g)_{<0>}
$$
for all $g\in H$, i.e. $\varphi:H\to A$ is left $H$-colinear,
where $A$ is a left $H$-comodule via $\tilde{\rho}$.\\
As opposed to the case ${}^H{\cal M}_A$, a left $H$-colinear map
$\varphi :H\to A$ is not sufficient to construct a quantum integral.
If however $\varphi :H\to A$ is a $k$-linear map
satisfying the more powerful relation
\begin{equation}\eqlabel{gal8}
\sum xg_{(1)} \ot \varphi(g_{(2)}) =
\sum S^{-1} \Bigl (\varphi (g)_{<1>}\Bigl)x
\varphi(g)_{<-1>}\ot \varphi(g)_{<0>}
\end{equation}
for all $x$, $g\in H$, then
$$
\gamma =\gamma_{\varphi}:H\to \Hom (H, A), \quad
\gamma (g)(h)=\varphi (S^{-1}(h)g)
$$
is a quantum integral. Indeed, the right hand side of \eqref{galTI1}
is
\begin{eqnarray*}
{\rm RHS}&=&
\sum S^{-1} \Bigl (\varphi ( S^{-1}(h_{(1)})g)_{<1>}\Bigl)
h_{(2)} \varphi ( S^{-1}(h_{(1)})g)_{<-1>}\ot
\varphi ( S^{-1}(h_{(1)})g)_{<0>}\\
\eqref{gal8}~~~&=&
\sum h_{(3)} S^{-1} (h_{(2)}) g_{(1)}\ot
\varphi ( S^{-1}(h_{(1)})g_{(2)})\\
&=& \sum g_{(1)}\ot \varphi ( S^{-1}(h) g_{(2)})\\
&=&\sum g_{(1)}\ot \gamma(g_{(2)})(h)
\end{eqnarray*}
i.e. $\gamma$ is a quantum integral. \\
2. The above remark indicates us a way to construct  quantum
integrals starting from integrals on $H$. Let $\theta :H\to k$
be a left integral on $H$. Then
$\varphi=\varphi_{\theta} :H\to A$, $\varphi (h)=\theta (h)1_A$
satisfies \eqref{gal8} and hence
\begin{equation}
\gamma=\gamma_{\theta}: H \to \Hom (H,A), \quad
\gamma (g)(h)=\theta (S^{-1}(h)g)1_A
\end{equation}
is a quantum integral. Furthermore, if $\theta (1_H)=1_k$ (that is
$H$ is cosemisimple), then $\gamma_{\theta}$ is a total quantum
integral.\\ 3. Assume that there exists $\gamma : H\to \Hom(H, A)$
a total quantum integral. Then any $M\in {}^H{\cal YD}_A$ is
relative injective as a left $H$-comodule. In particular, let $(A,
\rho^l, \rho^r)=(H,\Delta, \Delta)$, where $H$ is a finite
dimensional Hopf algebra over a field $k$. Then the extension $H^*
\subset D(H)$ is a RIT-extension, that is any representation of
the Drinfel'd double is injective as a right $H^*$-module.
\end{remarks}

\begin{proposition}\prlabel{mm4}
Let $H$ be a Hopf algebra with a bijective antipode and $A$ an
$H$-bicomodule algebra. Assume that there exists
$\gamma : H\to \Hom(H, A)$ a total quantum integral.
Then $\tilde{\rho}: A\to H\ot A$ splits in ${}^H{\cal YD}_A$.
\end{proposition}

\begin{proof}
Using \thref{doi} for ${}^H{\cal M}(H\ot H^{{\rm op}})_A={}^H{\cal YD}_A$,
the map
$$
\lambda :H\ot A\to A, \quad
\lambda (h\ot a)= \sum a_{<0>}\gamma (h)(S^{-1}(a_{<1>})a_{<-1>})
$$
for all $h\in H$, $a\in A$ is a left $H$-colinear retraction
of $\tilde{\rho}$. In particular, $\lambda (1_H\ot 1_A)=1_A$ and
\begin{equation}\eqlabel{galcol}
\sum g_{(1)}\ot \lambda (g_{(2)}\ot a)=
\sum S^{-1}\Bigl(\lambda (g\ot a)_{<1>}\Bigl)\lambda (g\ot a)_{<-1>}\ot
\lambda (g\ot a)_{<0>}
\end{equation}
for all $g\in H$ and $a\in A$. We define now
\begin{equation}\eqlabel{mm20}
\Lambda: H\ot A\to A, \quad
\Lambda (h\ot a)=\sum
\lambda\Bigl(S^{-2}(a_{<1>})hS(a_{<-1>})\ot 1_A\Bigl)a_{<0>}
\end{equation}
for all $h\in H$, $a\in A$. Then, for $a\in A$ we have
\begin{eqnarray*}
(\Lambda\circ\tilde{\rho})(a)&=&
\sum \Lambda (S^{-1}(a_{<1>})a_{<-1>}\ot a_{<0>})\\
&=&\sum \lambda \Bigl(S^{-2}(a_{<1>})S^{-1}(a_{<2>})a_{<-2>}S(a_{<-1>})\ot
1_A \Bigl)a_{<0>}\\
&=&\sum \lambda \Bigl(S^{-1}(a_{<2>}S^{-1}(a_{<1>}))a_{<-2>}S(a_{<-1>})\ot
1_A \Bigl)a_{<0>}\\
&=&\lambda (1_H\ot 1_A)a=a
\end{eqnarray*}
i.e. $\Lambda$ is still a retraction of $\tilde{\rho}$. Now,
for $h\in H$, $a$, $b\in A$ we have
\begin{eqnarray*}
\Lambda ((h\ot a)b)&=&\sum
\Lambda \Bigl(S^{-1}(b_{<1>})hb_{<-1>}\ot a b_{<0>} \Bigl)\\
&=&\sum
\lambda \Bigl(S^{-2}(a_{<1>})S^{-2}(b_{<1>})S^{-1}(b_{<2>})
hb_{<-2>}S(b_{<-1>})S(a_{<-1>}) \ot 1_A \Bigl)
a_{<0>}b_{<0>}\\
&=&\sum
\lambda \Bigl(S^{-2}(a_{<1>})S^{-1}(b_{<2>}S^{-1}(b_{<1>}))
hb_{<-2>}S(b_{<-1>})S(a_{<-1>}) \ot 1_A \Bigl)
a_{<0>}b_{<0>}\\
&=&\sum
\lambda \Bigl(S^{-2}(a_{<1>})
hS(a_{<-1>}) \ot 1_A \Bigl)
a_{<0>}b\\
&=& \Lambda (h\ot a)b
\end{eqnarray*}
hence $\Lambda$ is right $A$-linear. It remains to
prove that $\Lambda$ is also left $H$-colinear:
\begin{eqnarray*}
\tilde{\rho}\Lambda (h\ot a)
&=&\sum \tilde{\rho} \Bigl(
\lambda (S^{-2}(a_{<1>})hS(a_{<-1>}) \ot 1_A)a_{<0>}\Bigl)\\
&=&\sum S^{-1}(a_{<1>}) S^{-1} \Bigl(
\lambda (S^{-2}(a_{<2>})hS(a_{<-2>}) \ot 1_A)_{<1>} \Bigl)\\
&&\lambda (S^{-2}(a_{<2>})hS(a_{<-2>}) \ot 1_A)_{<-1>} a_{<-1>}\ot \\
&&\lambda (S^{-2}(a_{<2>})hS(a_{<-2>}) \ot 1_A)_{<0>} a_{<0>}\\
\eqref{galcol}~~
&=&\sum
S^{-1}(a_{<1>})S^{-2}(a_{<2>})h_{(1)}S(a_{<-2>}) a_{<-1>}\ot \\
&& \lambda (S^{-2}(a_{<3>})h_{(2)}S(a_{<-3>}) \ot 1_A)a_{<0>}\\
&=&\sum
S^{-1}\Bigl(S^{-1}(a_{<2>})a_{<1>} \Bigl)h_{(1)}S(a_{<-2>}) a_{<-1>}\ot \\
&& \lambda (S^{-2}(a_{<3>})h_{(2)}S(a_{<-3>}) \ot 1_A)a_{<0>}\\
&=& \sum
h_{(1)} \ot \lambda (S^{-2}(a_{<1>})h_{(2)}S(a_{<-1>}) \ot 1_A)a_{<0>}\\
&=& (Id\ot \Lambda)\rho_{H\ot A}(h\ot a)
\end{eqnarray*}
i.e. we proved that $\Lambda$ is a retraction of $\tilde{\rho}$
in ${}^H{\cal YD}_A$.
\end{proof}

We can define now the coinvariants of $A$ as follows:
$$
B=A^{{\rm co}(H)}:=\{\; a\in A \;\mid \; \tilde{\rho}(a)=1_H\ot a \;\}=
\{\; a\in A \;\mid \; \sum S^{-1}(a_{<1>}) a_{<-1>}\ot a_{<0>}=1_H\ot a \;\}
$$
Then $B$ is a subalgebra of $A$ and will be called the subalgebra of
{\em quantum coinvariants}.

\begin{proposition}\prlabel{mm5}
Let $H$ be a Hopf algebra with a bijective antipode, $A$ an
$H$-bicomodule algebra and $B$ the subalgebra of quantum
coinvariants. Assume that there exists $\gamma :H\to \Hom(H, A)$ a
total quantum integral.
Then:
\begin{enumerate}
\item $B$ is a direct summand of $A$ as a left $B$-submodule;
\item $B$ is a direct summand of $A$ as a right $B$-submodule.
\end{enumerate}
\end{proposition}

\begin{proof}
1. We shall prove that there exists a well defined left trace
given by the formula
\begin{equation}\eqlabel{mm200}
t^l: A\to B, \quad t^l (a)=\lambda (1_H\ot a)=
\sum a_{<0>}\gamma (1_H)(S^{-1}(a_{<1>})a_{<-1>})
\end{equation}
for all $a\in A$. Taking $g=1_H$ in \eqref{galcol} we obtain
$1_H \ot t^l (a)=\tilde{\rho}(t^l(a))$, i.e.
$t^l (a)\in B$, for all $a\in A$. Now, for $b\in B$
and $a\in A$
\begin{eqnarray*}
t^l (ba) &=& \sum b_{<0>}a_{<0>}\gamma (1_H) (S^{-1}(a_{<1>})
S^{-1}(b_{<1>}) b_{<-1>}a_{<-1>})\\ (b\in B)~~ &=& \sum
ba_{<0>}\gamma (1_H)(S^{-1}(a_{<1>})a_{<-1>})\\ &=& bt^l(a)
\end{eqnarray*}
hence $t^l$ is a left $B$-module map and finally
$$
t^l (1_A)=1_A\gamma (1_H)(1_H)=1_A\varepsilon (1_H)=1_A
$$
hence $t^l$ is a left $B$-module retraction of the inclusion
$B\subset A$.\\
2. Similarly, we can prove that the map given by the formula
\begin{equation}\eqlabel{mm21}
t^r: A\to B, \quad t^r (a)=\Lambda (1_H\ot a)=
\sum \gamma (S^{-2}(a_{<1>}) S(a_{<-1>}))(1_H) a_{<0>}
\end{equation}
for all $a\in A$, is a well defined right trace of the
inclusion $B\subset A$.
\end{proof}

\begin{definition}
Let $H$ be a Hopf algebra with a bijective antipode, $A$ an
$H$-bicomodule algebra and $\gamma :H\to \Hom(H, A)$ a
total quantum integral. The map
$$
t^l: A\to B, \quad t^l (a)=
\sum a_{<0>}\gamma (1_H)(S^{-1}(a_{<1>})a_{<-1>})
$$
for all $a\in A$ is called the (left) quantum trace associated to
$\gamma$.
\end{definition}

Now, we shall construct functors connecting ${}^H{\cal YD}_A$
and ${\cal M}_B$. First, if $M\in {}^H{\cal YD}_A$, then
$$
M^{{\rm co}(H)} =\{\; m\in M \;\mid \; \rho_M (m)=1_H\ot m \;\}
$$
is the right $B$-module of the coinvariants of $M$.
Furthermore, $M\to M^{{\rm co}(H)}$ gives us a covariant functor
$$
(-)^{{\rm co}(H)} : {}^H{\cal YD}_A\to {\cal M}_B.
$$
Now, for $N\in {\cal M}_B$, $N\ot_B A\in {}^H{\cal YD}_A$
via the structures
\begin{equation}\eqlabel{mm22}
(n\ot_B a)a'= n\ot_B aa'
\end{equation}
\begin{equation}\eqlabel{mm23}
\rho_{N\ot_B A}(n\ot_B a)=\sum S^{-1}(a_{<1>}) a_{<-1>}\ot n\ot_B a_{<0>}
\end{equation}
for all $n\in N$, $a$, $a'\in A$. In this way, we have constructed
a covariant functor called the induction functor
$$
-\ot_B A :{\cal M}_B \to {}^H{\cal YD}_A.
$$
We shall prove now that the above functors are an adjoint pair.
In the case $A=H$ the next result is the right-left version of the
Proposition 3.5 of \cite{CZ}.

\begin{proposition}
Let $H$ be a Hopf algebra with a bijective antipode and $A$ an
$H$-bicomodule algebra. Then
the induction functor $-\ot_B A :{\cal M}_B \to {}^H{\cal YD}_A$
is a left adjoint of the coinvariant functor
$(-)^{{\rm co}(H)} :{}^H{\cal YD}_A\to {\cal M}_B $.
\end{proposition}

\begin{proof}
Straightforward: the unit and the counit of the adjointness are
given by
\begin{equation}\eqlabel{mm24}
\eta_N :N \to (N\ot_B A)^{{\rm co}(H)}, \quad
\eta_N (n)=n\ot_B 1_A
\end{equation}
for all $N\in {\cal M}_B$, $n\in N$ and
\begin{equation}\eqlabel{mm25}
\beta_M :M^{{\rm co}(H)}\ot_B A \to M, \quad
\beta_M (m\ot_B a)=ma
\end{equation}
for all $M\in {}^H{\cal YD}_A$, $m\in M^{{\rm co}(H)}$ and $a\in A$.
\end{proof}

With the structures given by \eqref{gal9} and \eqref{gal101},
$H\ot A \in {}^H{\cal YD}_A$ and we identify $(H\ot A)^{{\rm co}(H)} \cong A$
via $ a\to 1_H\ot a$. Then the adjunction map $\beta_{H\ot A}$ can
be viewed as a map in ${}^H{\cal YD}_A$, as follows
\begin{equation}\eqlabel{mm229}
\beta=\beta_{H\ot A} :A\ot_B A\to H\ot A, \quad
\beta (a\ot_B b)=\sum S^{-1}(b_{<1>}) b_{<-1>}\ot ab_{<0>}
\end{equation}
for all $a$, $b\in A$. Here $A\ot_B A \in {}^H{\cal YD}_A$
via the structures
$$
(a\ot_B b)a'= a\ot_B ba', \quad
a\ot_B b \to \sum S^{-1}(b_{<1>}) b_{<-1>}\ot a\ot_B b_{<0>}
$$
for all $a$, $a'$, $b\in A$.

\begin{definition}
Let $H$ be a Hopf algebra with a bijective antipode, $A$ an
$H$-bicomodule algebra and $B=A^{{\rm co}(H)}$. Then $A/B$ is called
a quantum Galois extension if the canonical map
$$
\beta: A\ot_B A\to H\ot A, \quad
\beta (a\ot_B b)=\sum S^{-1}(b_{<1>}) b_{<-1>}\ot ab_{<0>}
$$
is bijective.
\end{definition}

If the right coaction $\rho^r :A\to A\ot H$ is trivial, the above
definition is just the left version of the usual well studied
$H$-Galois extensions (\cite{KT}, \cite{M}).
The quantum Galois extensions are the concepts that occur in the following
imprimitivity statement for quantum Yetter-Drinfel'd modules.

\begin{proposition}\prlabel{gal70}
Let $H$ be a Hopf algebra with a bijective antipode, $A$ an
$H$-bicomodule algebra and $B=A^{{\rm co}(H)}$. The following statements
are equivalent:
\begin{enumerate}
\item the induction functor
$-\ot_B A :{\cal M}_B \to {}^H{\cal YD}_A$
is an equivalence of categories;
\item the following conditions hold:
\begin{enumerate}
\item $A$ is faithfully flat as a left $B$-module;
\item $A/B$ is a quantum Galois extension.
\end{enumerate}
\end{enumerate}
\end{proposition}

\begin{proof}
$1)\Rightarrow 2)$ Trivial.\\ $2)\Rightarrow 1)$ is standard from
the categorical point of view: a pair of adjoint functors (as
$-\ot_B A :{\cal M}_B \to {}^H{\cal YD}_A$ and $(-)^{{\rm co}(H)}
:{}^H{\cal YD}_A\to {\cal M}_B $ are) gives an equivalence of
categories iff one of them is faithfully exact (or both of them
are exact) and the adjunction maps in the key objects of
categories ($B$ in ${\cal M}_B$ and $H\ot A$ in ${}^H{\cal YD}_A$)
are bijective. We point out that $\eta_N$ for all $N\in {\cal
M}_B$, and $\beta_M$ for all $M\in {}^H{\cal YD}_A$, can be
constructed from $\eta_B$ and $\beta_{H\ot A}$ using the
naturality condition: for details, in a more general frame, we
refer to Theorem 4.9 of \cite{CIMZ2} or \cite{CR} (for Doi-Koppinen
modules) or Theorem 3.10 of \cite{Br} (for entwining modules).
\end{proof}

We are going to prove now an affineness condition for quantum
Yetter-Drinfel'd
modules. First we need the following

\begin{theorem}\thlabel{mm70}
Let $H$ be a Hopf algebra with a bijective antipode, $A$ an
$H$-bicomodule algebra and $B=A^{{\rm co}(H)}$. Assume that there
exists $\gamma :H\to \Hom(H, A)$ a
total quantum integral. Then
$$\eta_N :N \to (N\ot_B A)^{{\rm co}(H)},\quad  \eta_N (n)=n\ot_B 1_A$$
is an isomorphism of right $B$-modules for all $N\in {\cal M}_B$.
\end{theorem}

\begin{proof}
Using the left quantum trace $t^l: A\to B$ we shall construct an
inverse of $\eta_N$. We define
\begin{equation}\eqlabel{mm260}
\theta_N :(N\ot_B A)^{{\rm co}(H)} \to N, \quad
\theta_N (\sum_i n_i\ot_B a_i)=\sum_i n_it^l(a_i)
\end{equation}
for all $\sum_i n_i\ot_B a_i\in (N\ot_B A)^{{\rm co}(H)}$.
As $t^l(1_A)=1_A$, $\theta_N\circ \eta_N=Id_N$. Let now,
$\sum_i n_i\ot_B a_i\in (N\ot_B A)^{{\rm co}(H)}$. Then
$$
\sum S^{-1}(a_{i_{<1>}}) a_{i_{<-1>}}\ot n_i\ot_B a_{i_{<0>}}=
\sum 1_H\ot n_i\ot_B a_i.
$$
It follows that (after we apply first $\gamma (1_H)$ and then
the flip map $\tau$)
$$
\sum n_i\ot_B a_{i_{<0>}} \ot
\gamma (1_H)(S^{-1}(a_{i_{<1>}}) a_{i_{<-1>}})=
\sum n_i\ot_B a_{i} \ot 1_A.
$$
Now, if we multiply the last factors we get
$$
\sum n_i\ot_B t^l (a_i)=\sum n_i\ot_B a_{i}.
$$
Hence we obtain
$$
(\eta_N\circ \theta_N)(\sum _i n_i\ot_B a_i)=
\sum_i n_it^l(a_i)\ot_B 1_A=
\sum n_i\ot_B t^l (a_i)=\sum n_i\ot_B a_{i}
$$
i.e. $\theta_N$ is an inverse of $\eta_N$.
\end{proof}

The above theorem applies in a large number of situations, as the
following Corollary shows.

\begin{corollary}
Let $A$ be an $H$-bicomodule algebra where $H$ is a
cosemisimple Hopf algebra over a field $k$ and $B=A^{{\rm co}(H)}$. Then
$$\eta_N :N \to (N\ot_B A)^{{\rm co}(H)},\quad  \eta_N (n)=n\ot_B 1_A$$
is an isomorphism of right $B$-modules for all $N\in {\cal M}_B$.
\end{corollary}

\begin{proof}
As $H$ is cosemisimple, there exists a left integral $\theta :H\to
k$ on $H$ with $\theta (1_H)=1_k$ (\cite{A}) and the antipode of
$H$ is bijective. Then, using 2) of \reref{galoi}, $\gamma :H\to
\Hom (H,A)$, $\gamma (g)(h)=\theta (S^{-1}(h)g)1_A$ for all $g$,
$h\in H$ is a total quantum integral and \thref{mm70} applies.
\end{proof}

\begin{remark}
\rm The above Corollary was proven recently in Theorem 3.4 of \cite{CZ}
in the case $(A=H, \rho^l=\rho^r=\Delta)$,
where $H$ is a semisimple and cosemisimple Hopf algebra.
The strategy adopted in \cite{CZ} for proving this result also
used the semisimplicity of $H$, which together with the cosemisimplicity
assures that $B=O(H)$ is a semisimple subalgebra of $H$.
\end{remark}

We shall prove now the main result of this section, that is the
affineness criterion for quantum Yetter-Drinfel'd modules.

\begin{theorem}\thlabel{mm80}
Let $H$ be a Hopf algebra with a bijective antipode and projective
over $k$, $A$ an $H$-bicomodule algebra and $B=A^{{\rm co}(H)}$.
Assume that:
\begin{enumerate}
\item there exists $\gamma : H\to \Hom(H, A)$ a total quantum
integral;
\item the canonical map
$\beta: A\ot_B A\to H\ot A, \;
\beta (a\ot_B b)=\sum S^{-1}(b_{<1>}) b_{<-1>}\ot ab_{<0>}$
is surjective.
\end{enumerate}
Then the induction functor
$-\ot_B A :{\cal M}_B \to {}^H{\cal YD}_A$
is an equivalence of categories.
\end{theorem}

\begin{proof}
In \thref{mm70} we have shown that, under the assumption 1),
the adjunction map
$\eta_N :N \to (N\ot_B A)^{{\rm co}(H)}$ is an
isomorphism for all $N\in {\cal M}_B$. It remains to prove that
the other adjunction map, namely
$\beta_M :M^{{\rm co}(H)}\ot_B A \to M$, $\beta_M (m\ot_B a)=ma$ is
an isomorphism for all $M\in {}^H{\cal YD}_A$. For this
we shall try to adapt the proof from Theorem 3.5 of
\cite{Sch}.\\
Let $V$ be a $k$-module. Then $A\ot V \in {}^H{\cal YD}_A$
via the structures induced by $A$ i.e.
\begin{eqnarray}
(a\ot v)b&=& ab\ot v \eqlabel{mm26}\\
\rho_{A\ot V}(a\ot v)&=& \sum S^{-1}(a_{<1>})a_{<-1>} \ot
a_{<0>}\ot v\eqlabel{mm27}
\end{eqnarray}
for all $a$, $b\in A$ and $v\in V$. In particular, for $V=A$,
$A\ot A \in {}^H{\cal YD}_A$ via
\begin{eqnarray}
(a\ot a')b&=& ab\ot a' \eqlabel{mm26'}\\
\rho_{A\ot A}(a\ot a')&=&
\sum S(a_{<1>}) a_{<-1>} \ot a_{<0>}\ot a'\eqlabel{mm27'}
\end{eqnarray}
for all $a$, $a'$, $b\in A$.
We will prove first that the adjunction map
$\beta_{A\ot V} :(A\ot V)^{{\rm co}(H)}\ot_B A \to A\ot V$ is an isomorphism
for any $k$-module $V$.\\
First, $V\ot B$ and $B\ot V\in {\cal M}_B$
via the usual $B$-actions $(v\ot b)b'=v\ot bb'$,
and  $(b\ot v)b'=bb'\ot v$ for all $v\in V$, $b$, $b'\in B$.
The flip map $\tau :V\ot B\to B\ot V$, $\tau (v\ot b)=b\ot v$ is an
isomorphism in ${\cal M}_B$.
On the other hand $V\ot A \in {}^H{\cal YD}_A$ via
\begin{eqnarray}
(v\ot a)b&=& v\ot ab \eqlabel{mm269}\\
\rho_{V\ot A}(v\ot a)&=&
\sum S^{-1}(a_{<1>}) a_{<-1>} \ot v\ot a_{<0>} \eqlabel{mm279}
\end{eqnarray}
for all $a$, $b\in A$ and $v\in V$. The flip map
$\tau :A\ot V\to V\ot A$, $\tau (a\ot v)=v\ot a$, is an isomorphism in
${}^H{\cal YD}_A$. Applying \thref{mm70} for $N=V\ot B\cong B\ot V$, we
obtain the following isomorphisms in ${\cal M}_B$
$$
B\ot V\cong V\ot B\cong (V\ot B\ot_B A)^{{\rm co}(H)}
\cong (V\ot A)^{{\rm co}(H)}\cong (A\ot V)^{{\rm co}(H)}
$$
The adjunction map $\beta_{A\ot V}$ for
$A\ot V$ is an isomorphism, as it is the composition of the
canonical isomorphisms
$$
(A\ot V)^{{\rm co}(H)}\ot_B A\cong (V\ot A)^{{\rm co}(H)}\ot_B A
\cong V\ot B\ot_B A\cong V\ot A \cong A\ot V.
$$
Let
$$
\tilde{\beta}: A\ot A\to H\ot A, \quad
\tilde{\beta} (a\ot b)=\sum S^{-1}(b_{<1>}) b_{<-1>}\ot ab_{<0>}
$$
for all $a$, $b\in A$. As $\beta$ is surjective, $\tilde{\beta}$
is surjective, because the diagram
\begin{diagram}
A\ot A    &                  &        &     &  \\
\dTo^{\rm can}      &  \SE^{\tilde{\beta}}   &        &      &  \\
A\ot_B A  & \rTo_{\beta}     & H\ot A & \rTo &0\\
\end{diagram}
is commutative, where ${\rm can}: A\ot A\to A\ot_B A$ is the canonical
surjection. Let us define now
\begin{equation}\eqlabel{mm28}
\zeta :A\ot A\to H\ot A, \quad
\zeta (a\ot b)=(\tilde{\beta}\circ \tau) (a\ot b)=
\sum S^{-1}(a_{<1>}) a_{<-1>}\ot ba_{<0>}
\end{equation}
for all $a$, $b\in A$. The map $\zeta$ is surjective, as
$\tilde{\beta}$ and $\tau$ are.
We will prove that $\zeta$ is a morphism
in ${}^H{\cal YD}_A$, where $A\ot A$ and $H\ot A$ are
quantum Yetter-Drinfel'd modules via \eqref{mm26'}, \eqref{mm27'} and
\eqref{gal9}, \eqref{gal101}.
Indeed,
\begin{eqnarray*}
\zeta ((a\ot b)a')&=& \zeta (aa'\ot b)\\
&=& \sum S^{-1}(a'_{<1>}) S^{-1}(a_{<1>}) a_{<-1>}a'_{<-1>}\ot
ba_{<0>}a'_{<0>}\\
\eqref{gal9}~~~
&=& \sum (S^{-1}(a_{<1>})a_{<-1>}\ot ba_{<0>})a'=\zeta (a\ot b)a'
\end{eqnarray*}
and
\begin{eqnarray*}
\rho_{H\ot A}(\zeta (a\ot b))&=&\sum
\rho_{H\ot A}(S^{-1}(a_{<1>})a_{<-1>}\ot ba_{<0>})\\
\eqref{gal101}~~~
&=&\sum S^{-1}(a_{<2>})a_{<-2>}\ot S^{-1}(a_{<1>}) a_{<-1>}\ot ba_{<0>}\\
&=&\sum (Id\ot \zeta)(S^{-1}(a_{<1>})a_{<-1>}\ot a_{<0>}\ot b)\\
&=& (Id\ot \zeta)\rho_{A\ot A}(a\ot b)
\end{eqnarray*}
for all $a$, $a'$, $b\in A$. Hence $\zeta$ is a surjective morphism
in ${}^H{\cal YD}_A$. \\
$H$ is projective over $k$; hence
$H\ot A$ is projective as a right $A$-module, where
$H\ot A$ is a right $A$-module in the usual way, $(h\ot a) b=h\ot ab$,
for all $h\in H$, $a$, $b\in A$. On the other hand, the map
$$
u: H\ot A\to H\ot A, \quad
u(h\ot a)= \sum S^{-1}(a_{<1>})h a_{<-1>}\ot a_{<0>}
$$
is an isomorphism of right $A$-modules: here the first $H\ot A$ has
the usual right $A$-module structure and the second $H\ot A$ has the
right $A$-module structure given by \eqref{gal9}. The $A$-linear
inverse of $u$ is given by
$$
u^{-1}: H\ot A\to H\ot A, \quad
u^{-1}(h\ot a)= \sum S^{-2}(a_{<1>})h S(a_{<-1>})\ot a_{<0>}
$$
In fact, $u$ is the isomorphism given by
\eqref{mm1.5}, associated to ${}^H{\cal YD}_A$.
We obtain that $H\ot A$, with the $A$-module structure
given by \eqref{gal9}, is still projective as a right $A$-module.
It follows that the surjective morphism $\zeta :A\ot A\to H\ot A$
splits in the category of right $A$-modules.
In particular, $\zeta$ is a $k$-split epimorphism in
${}^H{\cal YD}_A$.\\
Let now $M\in {}^H{\cal YD}_A$. Then
$A\ot A\ot M\in {}^H{\cal YD}_A$ via the structures arising from
$A\ot A$, that is
\begin{eqnarray}
(a\ot b\ot m)a'&=& aa'\ot b \ot m \eqlabel{mm261}\\
\rho_{A\ot A\ot M}(a\ot b\ot m)&=&
\sum S^{-1}(a_{<1>}) a_{<-1>} \ot a_{<0>}\ot b\ot m\eqlabel{mm271}
\end{eqnarray}
for all $a$, $b$, $a'\in A$, $m\in M$. On the other hand,
$H\ot A\ot M\in {}^H{\cal YD}_A$ via the structures arising
from the ones of $H\ot A$, i.e
\begin{eqnarray}
(h\ot a\ot m)b&=&
\sum S^{-1}(b_{<1>})hb_{<-1>}\ot ab_{<0>}\ot m\eqlabel{mm1400}\\
\rho_{H\ot A\ot M}(h\ot a\ot m)&=&
\sum h_{(1)}\ot h_{(2)}\ot a\ot m\eqlabel{mm14100}
\end{eqnarray}
for all $h\in H$, $a$, $b\in A$, $m\in M$. We obtain that $$ \zeta
\ot Id :A\ot A\ot M\to H\ot A\ot M $$ is a $k$-split epimorphism
in  ${}^H{\cal YD}_A$.\\ Applying \thref{mm3} for ${}^H{\cal
M}(H\ot H^{{\rm op}})_A={}^H{\cal YD}_A$ we obtain that the map $$
f: H\ot A\ot M\to M, \quad f(h\ot a\ot m)=\sum m_{<0>}\gamma
(S^{-2}(a_{<1>})h S(a_{<-1>}))(m_{<-1>})a_{<0>} $$ is a $k$-split
epimorphism in ${}^H{\cal YD}_A$. Hence, the composition $$
g=f\circ (\zeta \ot Id): A\ot A\ot M\to M,  \quad g(a\ot b\ot m)=
\sum m_{<0>}\gamma (S^{-2}(b_{<1>})S(b_{<-1>}))(m_{<-1>})b_{<0>}a
$$ is a $k$-split epimorphism in ${}^H{\cal YD}_A$. We note that
the structure of $A\ot A\ot M$ as an object in ${}^H{\cal YD}_A$
is of the form $A\ot V$, for the $k$-module $V=A\ot M$.\\ To
conclude, we have constructed a $k$-split epimorphism in
${}^H{\cal YD}_A$ $$ A\ot A\ot M=M_1 \stackrel{g}
{\longrightarrow}M{\longrightarrow}0 $$ such that the adjunction
map $\beta_{M_1}$ for $M_1$ is bijective. As $g$ is $k$-split and
there exists a total quantum integral $\gamma :H \to \Hom (H, A)$,
we obtain that $g$ also splits in ${}^H{\cal M}$. In particular,
the sequence $$ M_{1}^{{\rm co}(H)} \stackrel{g^{{\rm co}(H)}}
{\longrightarrow}M^{{\rm co}(H)}{\longrightarrow}0 $$ is exact.
Continuing the resolution with Ker($g$) instead of $M$, we obtain
an exact sequence in ${}^H{\cal YD}_A$ $$
M_2{\longrightarrow}M_1{\longrightarrow}M{\longrightarrow}0 $$
which splits in ${}^H{\cal M}$ and the adjunction maps for $M_1$
and $M_2$ are bijective. Using the Five lemma we obtain that the
adjunction map for $M$ is bijective.
\end{proof}

\thref{mm80} is a quantum affineness criterion.
We shall have a more complete picture of \thref{mmsch} for
quantum Yetter-Drinfel'd modules after solving the following open problem:

\vspace{5mm}
{\em Let $H$ be a Hopf algebra over a field $k$ with a bijective
antipode, $A$ an $H$-bicomodule algebra and $B=A^{{\rm co}(H)}$.
Assume that the induction functor
$-\ot_B A :{\cal M}_B \to {}^H{\cal YD}_A$
is an equivalence of categories. Is there a total quantum integral
$\gamma : H\to \Hom(H, A)$ ?}

\section{Conclusions and outlooks}
We have introduced the integrals associated to a Doi-Koppinen
datum $(H,A,C)$ and, as a major application, the quantum
integrals associated to the category of quantum Yetter-Drinfel'd modules.
The integrals of a Doi-Koppinen datum $(H, A, C)$
can be extended for an entwining structure
(\cite{BM}, \cite{Br}) or a weak entwining structure (\cite{CG}),
for a Doi-Koppinen datum over a weak
Hopf algebra (\cite{Bohm}) or for a Doi-Koppinen
datum $(H, A, C)$ in a braided category (\cite{BKLT}).

The notion of quantum Galois extensions was introduced by
quantizing the classic Galois extensions (\cite{KT}). The latter,
in the equivalent left version, can be retrieved by trivializing
the right coaction $\rho^r$. This process opens the way for
quantizing the Clifford theory of representations (\cite{Sch1})
and the theory of classic crossed products (\cite{BCM}). The
quantum Galois extensions introduced in the present paper are
related to quantum Yetter-Drinfel'd modules, while the quantum
Galois theory appearing in \cite{DM}, \cite{HMT} refers to vertex
operator algebras.

\vspace{5mm}
{\sl Acknowledgment.} All diagrams were drawn using the
"diagrams" software of Paul Taylor.

\end{document}